\documentclass[a4paper]{article}

\usepackage[latin1]{inputenc}
\usepackage{delarray,verbatim,enumerate}
\usepackage{amsmath,amsthm,amstext,amsbsy,amssymb,amsfonts,amscd}
\usepackage[frenchb]{babel}

\usepackage{ifpdf}
\ifpdf
	\usepackage{aeguill}
\else
	\usepackage[T1]{fontenc}
\fi

\newtheorem{Th}{Th\'eor\`eme}[section]
\newtheorem{Lem}[Th]{Lemme}
\newtheorem{Prop}[Th]{Proposition}
\newtheorem{Cor}[Th]{Corollaire}

\newtheorem{Sco}[Th]{Scolie}
\newtheorem{Def} [Th]{D\'efinition}

\def\Preuve{\smallskip\noindent {\it Preuve.~}}

\def\Remarque{\smallskip\noindent {\it Remarque.~}}

\font\teneufm=eufm10
\font\seveneufm=eufm7
\font\fiveeufm=eufm5
\newfam\gothfam
\textfont\gothfam=\teneufm \scriptfont\gothfam=\seveneufm
\scriptscriptfont\gothfam=\fiveeufm
\def\goth{\fam\gothfam}

\def\l{\ell}			\def\Tl{\mathbb  {T}_\ell}	\def\N{\mathbb N}
\def\Z{\mathbb {Z}}     			\def\Q{\mathbb Q}
\def\Zl{\mathbb{Z}_\ell}	\def\Fl{\mathbb{F}_\ell}	\def\Ql{\mathbb {Q}_\ell}	

\def\L{\Lambda}		\def\D{\Delta}			\def\g{\gamma}

\def\d{\partial}	 	
\def\J{\mathcal  J}  		\def\C{\mathcal  C}		\def\R{\mathcal  R}
\def\E{\mathcal  E} 		\def\G{\mathcal  G}		\def\U{\mathcal  U}
		
 			  	\def\Cl{\mathcal  C\l}

	\def\m{{\goth m}}		
	\def\p{{\goth p}}		

\def\p{\mathfrak p}		
\def\EE{{\mathfrak E}}	\def\PR{{\mathfrak R}}	\def\m{{\mathfrak m}}

\def\wi{\widetilde}		\def\ov{\overline}
	\def\deg{\operatorname{deg}}
\def\Gal{\operatorname{Gal}}	\def\Rad{\operatorname{Rad}}
	
\def\Hom{\operatorname{Hom}}

\title{\bf Généralisation d'un Théorème d'Iwasawa\footnote{J. Théor. Nombres Bordeaux {\bf 17} (2005), 527--553.}}

\author{Jean-François {\sc Jaulent}}

\date{}

\begin{document}

\maketitle

\centerline{\it A Georges Gras, à l'occasion de son soixantième anniversaire}\bigskip\bigskip

{\footnotesize
\noindent{\bf Résumé.} Nous généralisons à certains quotients finis d'un $\L$-module noethérien non nécessairement de torsion le classique théorème d'Iwasawa sur  l'expression asymptotique du $\l$-nombre de classes dans les $\Zl$-extensions. Puis nous illustrons les résultats obtenus en déterminant explicitement les caractères invariants attachés aux $\l$-groupes de $S$-classes $T$-infinitésimales dans une tour cyclotomique à partir de quelques paramètres référents et de données galoisiennes simples des extensions considérées. Un outil fondamental de cette étude est l'identité du miroir établie par Georges Gras, qui permet par dualité d'exprimer des conditions de ramification (sauvages ou modérées) dans une extension en termes de décomposition dans une autre extension. Les résultats obtenus précisent et complètent ceux établis dans un travail antérieur en collaboration avec Christian Maire (cf. \cite {jm}).
L'étude approfondie des contributions sauvages repose sur une généralisation d'un résultat de R. Greenberg (cf. \cite {gb}).}

\

{\footnotesize
\noindent{\bf Abstract.} We extend to convenient finite quotients of a noetherian $\L$-module the classical result of K. Iwasawa giving the asymptotic expression of the $\l$-part of the number of ideal class groups in $\Zl$-extensions of number fields. Then, in the arithmetic context, we compute the three characters associated by this way to the $\l$-groups of $T$-infinitesimal $S$-classes in the cyclotomic tower and relate them to the classical invariants and the decomposition characters associated to the finite sets of places $S$ and $T$. A main tool in this study is the so-called Spiegelungssatz of Georges Gras, which exchanges (wild or tame) ramification and decomposition. The main results of this arithmetical part extend those we obtained with Christian Maire in a previous article (cf. \cite {jm}). The most intricate study of the wild contribution of the sets $S$ and $T$ involves a generalization of a classical result of R. Greenberg on the genus theory of cyclotomic towers (cf. \cite {gb}).}

\


\noindent {\large \bf Introduction}\bigskip

Le résultat emblématique de la théorie d'Iwasawa (cf. \cite {Iw1} et \cite{Iw2})
affirme que les ordres respectifs $\l ^{x(n)}$ des $\l$-groupes de classes d'idéaux $\Cl_{K_n}$ attachés aux divers étages $K_n$ de la $\Zl$-extension cyclotomique $K_\infty = \cup_{n \in \N} K_n$ d'un corps de nombres $K$ sont donnés pour $n$ assez grand par une formule explicite de la forme~:
$$
x(n) \ =\ \mu \l ^n \ + \ \lambda n \ + \ \nu, 
$$
où $\nu$ est un entier relatif (éventuellement négatif) mais où $\lambda$ et $\mu$ sont des entiers naturels déterminés par la pseudo-décomposition de leur limite projective pour les applications normes $\,\C_{K_\infty} =  \varprojlim  \, \Cl_{K_n}$, regardée comme module de torsion sur l'algèbre d'Iwasawa $\Lambda = \Zl [[\gamma -1]]$ construite sur un générateur topologique $\gamma$ du groupe procyclique $\Gamma = \Gal(K_\infty/K)$ (cf. \cite {se}). Il est alors commode de reécrire l'égalité précédente sous une forme ne faisant intervenir que ces deux derniers paramètres~:
$$
x(n) \ \approx \ \mu \l ^n \ + \ \lambda n, 
$$
en convenant de tenir pour équivalentes deux suites d'entiers dont la différence est ultimement constante. L'identité obtenue vaut alors identiquement si l'on remplace les $\l$-groupes $\,\Cl_{K_n}$ par leurs quotients respectifs d'exposant $\l ^{n+1}$, comme expliqué dans \cite{ja0}.\medskip

L'objet de ce travail est précisément d'étendre ces résultats à des groupes de classes éventuellement infinis et dont la limite projective (pour la norme arithmétique) n'est donc plus généralement un $\L$-module de torsion.\medskip

Donnons nous pour cela deux ensembles finis disjoints $S$ et $T$ de places finies de $K$~; notons $S_n$ et $T_n$ respectivement les ensembles de places de $K_n$ au-dessus de $S$ et de $T$~; et considérons le $\l$-groupe $\Cl ^{T_n}_{S_n}$ des $S$-classes $T$-infinitésimales du corps $K_n$, que la théorie $\l$-adique du corps de classes identifie au groupe de Galois de la pro-$\l$-extension abélienne maximale $H^{T_n}_{S_n}$ de $K_n$ qui est $S$-décomposée et $T$-ramifiée (i.e. complètement décomposée aux places de $S_n$ et non ramifiée en dehors de $T_n$). Les groupes $\Cl ^{T_n}_{S_n}$ ne sont plus généralement finis (du moins dès que l'ensemble $T$ contient suffisamment de places au-dessus de $\l$), mais leurs quotients respectifs d'exposant $\l^{n+1}$ le sont encore, et il est montré dans \cite{jm} qu'il existe des entiers naturel $\rho^T_S$, $\mu^T_S$, $\lambda^T_S$ tels que les ordres respectifs $\l^{x^T_S(n)}$ des $^{\l^{n+1}}\!\Cl ^{T_n}_{S_n}$ vérifient l'estimation asymptotique~:
$$
x^T_S(n) \ \approx \ \rho^T_S (n+1)\l^n \ + \ \mu^T_S \l ^n \ + \ \lambda^T_S n.
$$
Ici encore les paramètres $\rho^T_S$, $\mu^T_S$ et $\lambda^T_S$ sont donnés par la pseudo-décomposi\-tion élémentaire de la limite projective $\C^{S_\infty}_{T_\infty} = \varprojlim \ \Cl^{T_n}_{S_n}$ regardée comme module noethérien sur l'algèbre d'Iwasawa $\Lambda$, à cela près que l'invariant $\lambda^T_S$ diffère d'une unité de l'invariant lambda du $\Lambda$-module $\C^{T_\infty}_{S_\infty}$ lorsque l'extension $K_{\infty}/K$ est elle-même $S$-décomposée et $T$-ramifiée, i.e. lorsque l'ensemble $S$ est vide et que $T$ contient la totalité des places au-dessus de $\l$.\medskip

De fait, la démonstration de la formule précédente repose intégralement sur des résultats de \cite{ja0} (essentiellement ceux du chapitre IV) qui n'avaient jamais fait l'objet jusqu'ici d'une publication spécifique, et même en toute rigueur sur une version plus forte non publiée de ces résultats prenant en compte les $\Lambda$-modules noethériens qui ne sont pas nécessairement de torsion.\par 

Il nous a donc paru utile d'en rassembler ici une preuve complète. Et c'est l'objet de la première partie, purement algébrique, de ce travail.\medskip

De plus, une fois la formule établie, se pose évidemment, dans le contexte arithmétique évoqué plus haut, la question du calcul des divers paramètres des $\l$-groupes de $S$-classes $T$-infinitésimales, ou du moins de leur détermination à partir de quelques {\it paramètres référents} ne faisant pas intervenir les ensembles de places $S$ et $T$ et d'invariants galoisiens simples qui leur sont attachés. Et c'est l'objet de la seconde partie, arithmétique, de ce travail.\smallskip

Or dans cette dernière étude, les {\it identités du miroir} obtenues par Georges Gras (cf. \cite {gr1}) se révèlent un outil essentiel. Et comme celles-ci s'expriment naturellement en termes de représentations, il est préférable, pour les utiliser pleinement, de remplacer les paramètres traditionnels d'Iwasawa, qui s'interprètent comme dimensions de certains espaces vectoriels, par des {\it caractères} prenant explicitement en compte l'action galoisienne et dont les degrés respectifs ne sont autres que les invariants de départ. Et c'est donc en termes de caractères que nous énonçons les théorèmes de paramétrage.\medskip

Une conséquence immédiate des identités de dualité est le déploiement (déjà remarqué dans \cite{jm}) de la ramification modérée au-dessus de la tour cyclotomique. En revanche, l'étude fine de la contribution des places sauvages soulève rapidement de délicats problèmes d'indépendance $\l$-adique. La formule finale que nous obtenons vaut ainsi sous la conjecture de Leopoldt. Elle repose sur une généralisation dans le cadre $\l$-ramifié d'un résultat de semi-simplicité de R. Greenberg (cf. \cite{gb}), dont nous donnons en appendice une démonstration autonome basée sur la Théorie des Genres.\smallskip 

En fin de compte, il résulte de cette étude que les divers paramètres d'Iwasawa respectivement associés aux groupes de Galois $\C^{T_\infty}_{S_\infty}$ attachés aux pro-$\l$-extensions abéliennes $T$-ramifiées $S$-décomposées maximales sur une tour cyclotomique sont donnés par des formules entièrement explicites à partir des seuls paramètres {\it réels} associés aux seuls groupes $\C^{T_\infty^\l} = \C_\emptyset^{T^\l}$, où $T^\l$ décrit les sous-ensembles de l'ensemble $L$ des places $\l$-adiques du corps de base considéré.
\bigskip


\noindent {\large \bf Principales notations}\bigskip

Nous utilisons dans ce qui suit les notations et conventions suivantes~:

\bigskip

\noindent{\it Notations algébriques de la Théorie d'Iwasawa}

\smallskip\begin{itemize}

\item$\l$ est un nombre premier arbitraire (pris impair dans la deuxième section)~;
\item$\mu_{\l^\infty} = \cup_{n \in \N}\, \mu_{\l^n}$ est le groupe des racines d'ordre $\l$-primaire de l'unité~;
\item$\D$ est un groupe abélien fini d'ordre $d$ étranger avec $\l$~;
\item$\Zl [\Delta ] \, =\, \oplus _\varphi \ \Zl [\Delta ] e_\varphi\, =\, \oplus _\varphi \ Z_\varphi$ est l'algèbre $\l$-adique du groupe $\D$~;
\item$\Gamma = \gamma^{\Z_\l}$ est le groupe procyclique isomorphe à $\Z_\l$ engendré par $\gamma$~;
\item$\Lambda = \Z_\l[[\gamma - 1]]$ est l'algèbre d'Iwasawa associée au groupe $\Gamma$~;
\item$\omega_n$ est le polynôme $\gamma^{\l^n}-1$ de l'anneau $\Zl[\gamma -1]=\Zl[\gamma]$~;
\item$\nabla_n$ est l'idéal de $\Lambda$ engendré par le polynôme $\omega_n$ et l'élément $\l^{n+1}$~;
\item$\Lambda[\Delta] = \oplus_\varphi \,\Lambda[\Delta] e_\varphi = \oplus_\varphi\, \Lambda_\varphi$ est l'algèbre du groupe $\D$ sur l'anneau $\L$.

\end{itemize}\bigskip

\noindent{\it Notations arithmétiques de la Théorie d'Iwasawa}

\smallskip\begin{itemize}

\item$F$ est un corps de nombres totalement réel~;
\item$F_\infty = \cup_{n \in \N}\, F_n$ est la $\Zl$-extension cyclotomique de $F$~;
\item$K$ est une extension abélienne de $F$ de groupe $\D$ et contenant $\mu_\l$~;
\item$\D = \Gal(K/F)$ est un groupe abélien fini d'ordre $d$ étranger avec $\l$~;
\item$K_{\infty} = \cup_{n \in \N}\, K_n$ est la $\Z_\l$-extension cyclotomique de $K$ (et contient $\mu_{\l^\infty}$)~;
\item$K_n$ est l'unique sous-corps de $K_\infty$ de degré $\l^n$ sur $K$ (on a ainsi $K = K_0$)~;
\item$\Gamma = \gamma^{\Z_\l}$ est le groupe de Galois $\Gal(F_{\infty}/F) \simeq \Gal(K_{\infty}/K)$~;
\item$\Gamma_n = \gamma^{\l^n\Z_\l}$ est le sous-groupe $\Gal(F_{\infty}/F_n) \simeq \Gal(K_{\infty}/K_n)$~;
\item$\Lambda = \Z_\l[[\gamma - 1]]$ est l'algèbre d'Iwasawa associée à $\Gamma$~;

\medskip

\item$S$ et $T$ sont deux $\D$-ensembles finis disjoints de places finies de $K$~;
\item$L$ est l'ensemble des places $\l$-adiques (i.e. au-dessus de $\l$) du corps $K$~;
\item$S^\l =S \cap L$ est la partie sauvage de $S$ et $S^0 = S \setminus S_\l$ sa partie modérée~;
\item$T^\l =T \cap L$ est la partie sauvage de $T$ et $T^0 = S \setminus T_\l$ sa partie modérée~;
\item$S_n$ et $T_n$ sont les ensembles finis de places de $K_n$ au-dessus de $S$ et de $T$~;
\item$S_\infty$ et $T_\infty$ sont les ensembles finis de places de $K_\infty$ au-dessus de $S$ et $T$~;

\item
$\deg _\l S = \sum_{\p \in S_{\l}}[F_{\p}:\Q_\l]$ est le degré sauvage en $S$ ; de même $\deg _\l T$ en $T$.

\end{itemize}\bigskip

\noindent{\it Notations de la Théorie $\l$-adique du Corps de Classes}

\smallskip\begin{itemize}

\item $\p$ est une place finie du corps $K$~et  $\p_n$ une place de $K_n$ au-dessus de $\p$~;
\item $K_{\p_n}$ : le complété de $K_n$ en la place $\p_n$~;
\item $\R_{\p_n} = \varprojlim K^{\times}_{\p_n}/K^{\times {\l}^m}_{\p_n}$ : le compactifié $\ell$-adique du groupe $K^{\times}_{\p_n}$~;
\item $\U_{{\goth p}_n}$ : le sous-groupe unité de $\R_{{\goth p}_n}$~; 
\item $\mu_{{\goth p}_n}$ : le sous-groupe de torsion de $\R_{\p_n}$, i.e. le groupe des racines de l'unité~;\medskip

\item $\J_{K_n} = \prod^{\rm res}_{{\p}_n} \R_{{\goth p}_n}$ : le $\l$-groupe des idèles du corps $K_n$~;
\item $\U_{K_n} = \prod_{{\p}_n} \U_{{\goth p}_n}$ : le sous-groupe des idèles unités~;
\item $\J^{T_n} = \prod_{\p_n \in T_n}^{\rm res} \R_{\p_n}$ : le groupe des idèles $T$-infinitésimaux~;
\item $\U^{T_n} =\, \prod_{\p_n \in T_n} \U_{\p_n}$ : le sous-groupe unité de $\J^{T_n}$~;
\item $\J_{S_n}^{T_n} =\, \prod_{\p_n \in S_n} \J_{\p_n} \prod_{\p_n \notin T_n} \U_{\p_n}$ : le groupe des $S$-idèles $T$-infinitésimaux~;\medskip

\item $\R_{K_n} = \Zl \otimes_{\Z} K^{\times}_n \subset \J_{K_n}$ : le sous-groupe principal de $\J_{K_n}$~;
\item $\E_{K_n}  = \R_{K_n}\! \cap \,\U_{K_n} = \Zl \otimes_{\Z} E_{K_n}$ : le $\l$-adifié du groupe des unités globales~;
\item $\E^{T_n}_{K_n}  = \R_{K_n}\! \cap \,\U^{T_n}_{K_n}$ : le sous-groupe des unités $T$-infinitésimales de $K_n$~;
\item $\E^{T_n}_{S_n}  = \R_{K_n}\! \cap \,\J^{T_n}_{S_n}$ : le groupe des $S$-unités $T$-infinitésimales de $K_n$~;\medskip

\item $\Cl_{K_n} = \J_{K_n}/ \U_{K_n} \R_{\p_n}$ : le $\l$-groupe des classes (au sens habituel) de $K_n$~;
\item $\Cl_{S_n}^{T_n} = \J_{K_n}/ \J_{S_n}^{T_n} \R_{\p_n}$ : le $\l$-groupe des $S$-classes $T$-infinitésimales de $K_n$~;
\item $\Cl_{S_\infty}^{T_\infty} = \varinjlim \,\Cl_{S_n}^{T_n}$ : le groupe des $S$-classes $T$-infinitésimales de $K_\infty$~;
\item $\C_{S_\infty}^{T_\infty} = \varprojlim \,\Cl_{S_n}^{T_n}$ : la limite projective (pour la norme) des groupes $\,\Cl_{S_n}^{T_n}$~;
\item $\rho^T_S, \mu^T_S, \lambda^T_S$ : les caractères invariants d'Iwasawa du $\Lambda[\D]$-module $\,\C_{S_\infty}^{T_\infty}$~;

\end{itemize}\bigskip

\noindent{\it Notations relatives à l'involution du miroir}

\medskip\begin{itemize}

\item $\mathbb T_\l =\varprojlim \mu_{\l^n}$ est le module de Tate construit sur les racines de l'unité~;
\item $\D_\p$ : le sous-groupe de décomposition des $\p_n$ dans $K_n/F_n$ pour $n\gg  0$~~ ;
\item $\chi_\p$ : l'induit à $\D$ du caractère unité de $\D_\p$~;
\item $\chi_T = \sum_{\p \in T} \chi_\p$, la somme étant prise sur les places de $F_\infty$~;
\item $\omega$ : le caractère cyclotomique, i.e. le caractère du $\Zl[\D]$-module $\mathbb T_\l$~;
\item $\psi \to \psi^* = \omega \psi^{-1}$ est l'involution du miroir (où $\psi^*$ est le {\it reflet} de $\psi$)~;
\item $\chi_{\rm reg}$ : le caractère régulier de $\D$~;
\item $\chi = \chi^\oplus + \chi^\ominus$ : la décomposition de $\chi$ en ses parties réelles et imaginaires~;
\item $\chi_\infty = [F : \Q] \,\chi_{\rm reg}^{ \, \oplus}$ et $\delta_S =\deg _\l S \, \chi_{\rm reg}$~; de même~: $\delta_T =\deg _\l T \, \chi_{\rm reg}$~; 

\medskip

\item $\m_n = \prod_{\p_n \in T} \p_n^{\nu_\p}$ : un diviseur convenable d'hyperprimarité du corps $K_n$~;
\item $\R^{\m_n}$ : le sous-module de $\R_{K_n}\!$ formé des éléments congrus à 1 modulo $\m_n$~;
\item $\E^{\m_n}_{S_n} = \E_{S_n} \cap \R^{\m_n}$ : le groupe des $S$-unités congrues à 1 modulo $\m_n$~;\smallskip

\item $\Cl^{\m_n}_{S_n} = \J_{K_n} / \J^{\m_n}_{S_n} \R_{K_n}$ le $\l$-groupe des $S$-classes de rayon $\m_n$ dans $K_n$.

\item $H^{T_n}_{S_n}$ : l'extension abélienne $S$-décomposée $T$-ramifiée d'exposant $\l^{n+1}$ maximale du corps $K_n$~;

\item $\Rad ^{T_n}_{S_n}$ : le radical kummérien de l'extension abélienne $H^{T_n}_{S_n}/ K_n$~;

\item $\Gal ^{T_n}_{S_n}$ : le groupe de Galois de l'extension abélienne $H^{T_n}_{S_n}/ K_n$.

\end{itemize}\
\newpage


\setcounter{section}{1}  \setcounter{Th}{0} 

\noindent{\large \bf 1. Le Théorème du paramétrage pour les $\L[\D]$-modules}\bigskip


\noindent{\bf 1.1 Rappel du contexte algébrique de la Théorie d'Iwasawa}\bigskip

Dans cette section purement algébrique, $\l$ est un nombre premier arbitraire, $\Zl$ désigne l'anneau local des entiers $\l$-adiques, $\Fl \simeq \Z / \l\Z$ son corps résiduel et $\Delta$ est un groupe abélien fini d'ordre $d$ étranger à $\l$.\smallskip

 Sous l'hypothèse $\ell \nmid d$, l'algèbre résiduelle $\Fl [\Delta ]$ est ainsi une algèbre semi-simple, produit direct d'extensions $F_\varphi $ de $\Fl$~; l'algèbre $\ell$-adique $\Zl [\Delta ]$ est une algèbre semi-locale, produit direct d'extensions non ramifiées $Z_\varphi $ de $\Zl$~; et les idempotents primitifs $\bar e_\varphi$  (respectivement $e_\varphi $) correspondant à leurs décompositions respectives~:
$$
\Fl [\Delta ]\, =\, \oplus _\varphi \ \Fl [\Delta ]\bar e_\varphi\, =\, \oplus _\varphi \ F_\varphi \qquad \& \qquad
\Zl [\Delta ] \, =\, \oplus _\varphi \ \Zl [\Delta ] e_\varphi\, =\, \oplus _\varphi \ Z_\varphi
$$
sont donnés à partir des caractères $\ell$-adiques irréductibles $\varphi $ de $\Delta $ par les formules classiques~:

\centerline{$e_\varphi\, =\, \frac{1}{d}\sum_{\tau\in \Delta}\varphi(\tau^{-1})\tau$,}\medskip

\noindent ainsi que leurs réductions respectives modulo $\l$. En particulier, la dimension $[\Z_\varphi :\Zl]$ du $\Zl$-module $\Z_\varphi$ (c'est à dire le degré de l'extension associée $\Q_\varphi / \Ql$ des corps de fractions) est le degré $\deg \varphi$ du caractère $\varphi$.\smallskip

De façon toute semblable, si $\Gamma = \g^{\Zl}$ désigne un groupe procyclique isomorphe à $\Zl$ et  $\Lambda = \Zl [[\gamma -1]]$ l'algèbre des séries formelles en l'indéterminée $\gamma -1$ à coefficients dans l'anneau $\Zl$ des entiers $\l$-adiques, l'algèbre de groupe $\Lambda[\Delta]$ s'écrit canoniquement comme produit direct de ses $\varphi$-composantes~:
$$
\Lambda[\Delta] = \oplus_\varphi \,\Lambda[\Delta] e_\varphi = \oplus_\varphi\, \Lambda_\varphi\ .
$$
Et, pour chaque caractère irréductible $\varphi$ du groupe $\Delta$, la $\varphi$-composante $\Lambda_\varphi = \Lambda[\Delta] e_\varphi$ associée à l'idempotent $e_\varphi$ s'identifie à l'algèbre des séries formelles $\Z_\varphi [[\gamma -1]]$ en l'indéterminée $\gamma -1$ à coefficients dans l'extension non ramifiée $\Z_\varphi =\Zl[\Delta] e_\varphi$ de l'anneau $\Zl$.\medskip

Plus généralement, par action des idempotents primitifs $e_\varphi$ de l'algèbre $\L[\D]$ tout $\L$-module noethérien $X$ se décompose naturellement comme somme directe de ses $\varphi$-composantes $X_\varphi = X^{e_\varphi}$, c'est à dire comme somme directe de $\L_\varphi$-modules noethériens. Les théorèmes de structure de la Théorie d'Iwasawa, tels qu'énoncés par Serre (cf. \cite{se}), permettent donc {\it mutatis mutandis} de décrire un tel module à pseudo-isomorphisme près à partir d'un certain nombre de modules référents dits élémentaires. Il vient ainsi~:

\begin{Th} Tout $\L[\D]$-module noethérien $X$ est pseudo-isomorphe à une somme directe finie de modules isotypiques élémentaires.\par

Plus précisément pour chaque caractère $\l$-adique irréductible $\varphi$ du groupe $\D$, il existe un unique triplet $(\rho_\varphi,\ s_\varphi,\ t_\varphi)$ d'entiers naturels, une unique suite décroissante $(f_{\varphi ,i})_{i=0,\dddot \ ,s_\varphi}$ de polynômes distingués de l'anneau $\Z_\varphi[\g-1]$ et une unique suite décroissante $(m_{\varphi ,i})_{i=0,\dddot \ ,t_\varphi}$d'entiers naturels non nuls tels que la $\varphi$-composante $X_\varphi = e_\varphi X$ du module $X$ s'envoie par un $\L[\D]$-morphisme à noyau et conoyau fini dans la somme directe~:
$$
X_\varphi \sim \L_\varphi^{\rho_\varphi} \oplus \big( \oplus_{i=0}^{s_\varphi}\L_\varphi /f_{\varphi,i}\L_\varphi \big) \oplus \big( \oplus_{j=0}^{t_\varphi}\L_\varphi / \l^{m_{\varphi,i}}\L_\varphi \big).
$$

On dit que l'entier $\rho_\varphi = \dim_{\L_\varphi}X_\varphi$ est la dimension du $\L[\D]$-module $X_\varphi$ et que le polynôme $P_\varphi = \prod_{j=0}^{t_\varphi}\l^{m_{\varphi,j}}\prod_{i=0}^{s_\varphi}f_{\varphi,i} \ \in \Z_\varphi[\g-1]$ est le polynôme caractéristique de son sous-module de $\L_\varphi$-torsion.
\end{Th}

Il est alors commode de coder globalement l'information  donnée par l'ensemble de ces invariants en introduisant comme suit la notion de paramètres~:

\begin{Def} Dans le contexte précédent, nous appelons paramètres d'un  $\L[\D]$-module noethérien $X$ les caractères $\l$-adiques du groupe abélien $\D$ définis à partir des invariants d'Iwasawa des composantes isotypiques de $X$ par les formules~:\smallskip

\centerline{$\rho = \sum_\varphi \rho_\varphi \ \varphi \ , \qquad
\mu = \sum_\varphi \mu_\varphi \ \varphi \ , \qquad
\lambda = \sum_\varphi \lambda_\varphi \ \varphi $,}\medskip

\noindent où, pour chaque caractère $\l$-adique irréductible $\varphi$, les entiers $\rho_\varphi$, $\mu_\varphi$ et $\lambda_\varphi$ mesurent respectivement la dimension $\dim_{\L_\varphi}X_\varphi$ de la $\varphi$-composante de $X$, ainsi que la $\l$-valuation $\sum_{j=0}^{t_\varphi}m_{\varphi,j}$ et le degré $\sum_{i=0}^{s_\varphi}\deg f_{\varphi,i}$ du polynôme caractéristique de son sous-module de $\L_\varphi$-torsion.
\end{Def}

\medskip
Introduisons maintenant la filtration $(\nabla_{\!n})_{n \in \N}$, qui joue un rôle essentiel dans ce qui suit~:

\begin{Prop} Pour chaque entier naturel $n$, désignons par \smallskip

\centerline{$\nabla_{\!n} = \l^{n+1} \L + \omega_n \L$}\smallskip

\noindent l'idéal de l'algèbre d'Iwasawa $\L = \Z_\l[[\g-1]]$ engendré par l'élément $\l^{n+1}$ et le polynôme $\omega_n =\g^{\l^n}-1$ de l'anneau $\Z_\l[\g-1] = \Z_\l[\g]$. Cela étant~:
\begin{enumerate} 

\item[(i)] L'idéal $\nabla = \nabla_{\!0}$ est l'unique idéal maximal de l'anneau local régulier $\L$.

\item[(ii)] Les idéaux $(\nabla_{\!n})_{n \in \N}$ forment une suite exhaustive strictement décroissante d'idéaux d'indice fini de l'algèbre d'Iwasawa $\L = \Z_\l[[\g-1]]$.

\item[(iii)] Un $\L[\D]$-module compact $X$ est noethérien si et seulement s'il existe un entier $n$ pour lequel le quotient $X / \nabla_{\!n} X$ est fini~; auquel cas, les quotients $X / \nabla_{\!n} X$ sont tous des $\L[\D]$-modules finis.

\item[(iv)]
Un $\L[\D]$-module compact $X$ est fini si et seulement s'il existe un entier $n$ pour lequel le sous-module $X / \nabla_{\!n} X$ est trivial~; auquel cas, les sous-modules $X / \nabla_{\!n} X$ sont ultimement triviaux.
\end{enumerate} 
\end{Prop}

\Preuve Les deux premières assertions s'obtiennent comme suit~: d'un côté, 
comme $\omega_0$ est égal à $\gamma -1$, l'idéal $\nabla_{\!0}$ est précisément 
l'idéal maximal $\nabla$ de l'algèbre $\Lambda$~; et l'identité\smallskip

\centerline{$(\gamma^{\l^n} - 1)^\l = (\gamma^{\l^{n+1}} - 1) + \l (\gamma^{\l^n} - 1) P(\gamma^{\l^n})$,}\medskip

\noindent pour un polynôme convenable $P$ de l'anneau $\Zl [\gamma -1] = \Zl [\gamma]$,
montre par une récurrence évidente que $\nabla_{\!n}$ est contenu dans la 
puissance $(n+1)$-ième de l'idéal $\nabla$. D'un autre côté, un calcul 
immédiat donne~:\medskip

\centerline{$(\Lambda : \nabla_{\!n}) = |(\Z / \l^{n+1}\Z)[\Gamma_n]| = \l^{(n+1)\l^n}$,}\medskip

\noindent avec $\Gamma_n = \Gamma/ \Gamma^{\l^n}$~; ce qui établit le deuxième point. \par
Les assertions suivantes résultent alors du fait que les idéaux $(\nabla_{\!n})_{n \in \N}$ forment une base de voisinages de 0 dans l'algèbre topologique compacte $\Lambda = \Zl [[\gamma -1]]$.\medskip

L'objet de la section qui suit est de relier les paramètres invariants d' un $\L[\D]$-module noethérien aux ordres de certains de ses quotients finis.

\bigskip

\noindent{\bf 1.2  Enoncé du théorème principal sur le paramétrage}\bigskip

Le résultat essentiel de cette section relie les paramètres d'un $\L[\D]$-module noethérien donné $X$ avec les ordres de certains sous-quotients finis de $X$~:\smallskip

\begin{Th} {\rm{\bf (Théorème des paramètres.)}} Soit $X$ un $\L[\D]$-module noethérien de  paramètres
$\rho$, $\mu$ et $\lambda$. Si $\nabla_{\!n} = \l^{n+1} \L + \omega_n \L$ désigne
l'idéal de l'algèbre d'Iwasawa $\L = \Z_\l[[\g-1]]$ engendré par l'élément
$\l^{n+1}$ et le polynôme $\omega_n =(\g^{\l^n}-1)$, il existe un unique
caractère $\l$-adique virtuel $\nu$ du groupe $\D$ tel que l'ordre $\l^{x_n^\varphi}$
de la $\varphi$-composante du quotient $X/ \nabla_{\!n} X$ soit donné, pour chaque
caractère $\l$-adique irréductible $\varphi$ et tout entier $n$ assez grand, par
la formule~:\smallskip

\centerline{$x_n^\varphi \ = \ <\rho, \varphi>(n+1)\l^n \ + \ <\mu, \varphi>\l^n \ + \ <\lambda, \varphi>n \ + \ <\nu, \varphi>$.}\smallskip

En d'autres termes, l'indice $\l^{x_n^\varphi} = \big( X_\varphi : \nabla_{\!n} X_\varphi \big)$ est donné asymptotiquement par la formule~:

\centerline{$x_n^\varphi \ \approx \ <\rho, \varphi>(n+1)\l^n \ + \ <\mu, \varphi>\l^n \ + \ 
<\lambda, \varphi>n$ ,}\smallskip

\noindent en ce sens que la différence entre les deux membres est ultimement constante.
\end{Th}\smallskip

Ce résultat amène ainsi à poser~:\smallskip

\begin{Def} Nous disons qu'une suite $(X_n)_{n \in \N}$ de $\Zl [\Delta]$-modules 
finis est paramétrée par les caractères $\l$-adiques virtuels~:\smallskip

\centerline{$\rho = \sum_\varphi \rho_\varphi \ \varphi \ , \qquad
\mu = \sum_\varphi \mu_\varphi \ \varphi \ , \qquad
\lambda = \sum_\varphi \lambda_\varphi \ \varphi \qquad 
\nu = \sum_\varphi \nu_\varphi \ \varphi $,}\smallskip

\noindent lorsque pour chaque caractère $\l$-adique irréductible $\varphi$ du groupe $\D$ l'ordre $\l^{x_n^\varphi}$ de la $\varphi$-composante $e_\varphi X_n$ de $X_n$ est donné asymptotiquement par la formule~:\smallskip

\centerline{$x_n^\varphi \ = \ <\rho, \varphi>(n+1)\l^n \ + \ <\mu, \varphi>\l^n \ + \ <\lambda, \varphi>n \ + \ <\nu, \varphi>$.}\smallskip

Lorsque le caractère $\nu$ n'est pas explicité, nous écrivons simplement~:\smallskip

\centerline{$x_n^\varphi \ \approx \ <\rho, \varphi>(n+1)\l^n \ + \ <\mu, \varphi>\l^n \ + \ 
<\lambda, \varphi>n$,}\smallskip

\noindent pour signifier que la différence entre les deux membres est ultimement constante.
\end{Def}

\Remarque L'entier $x_n^\varphi$ étant toujours positif, le caractère $\rho$ qui apparait dans la Définition 1.5 est de ce fait positif (en ce sens que l'on a $\rho_\varphi \ge 0$ pour chaque caractère irrédutible $\varphi$) ; il en est de même du caractère $\mu$ si $\rho$ est nul ; et du caractère $\lambda$ si $\rho$ et $\mu$ sont tous deux nuls. Il peut arriver en revanche que les caractères $\lambda$ ou $\mu$ ne soient pas dans $R^+_{\Zl} (\Delta)$ dès lors que le caractère dominant y est~: des exemples sont donnés dans la section 2 plus loin.\medskip

\noindent {\it Preuve du Théorème dans le cas élémentaire.} Distinguons les trois cas~:\smallskip
\begin{itemize}
\item[(i)] Pour $X=\L_\varphi$, il vient tout simplement~:\smallskip

\centerline{$X/ \nabla_{\!n} X \simeq (\Z_\varphi / \l^{n+1}\Z_\varphi)[\Gamma_n]$,  avec $\Gamma_n = \Gamma / \Gamma^{\l^n} \simeq \Z/ \l^n \Z$~;}\smallskip

\noindent d'où l'identité attendue~:\smallskip

\centerline{$(X : \nabla_{\!n} X) = (\Z_\varphi : \l^{n+1}\Z_\varphi)^{\l^n} = \l^{(n+1) \l^n \deg\varphi}$.}\smallskip

\item[(ii)] Pour $X = \L_\varphi / \l^{m_\varphi}\L_\varphi$ et $n \ge m_\varphi$, le même calcul donne~:\smallskip

\centerline{$(X : \nabla_{\!n} X) = (\Z_\varphi : \l^{m_\varphi}\Z_\varphi)^{\l^n} = \l^{m_\varphi \l^n \deg\varphi}$.}\smallskip

\item[(iii)] pour $X = \L_\varphi / f_\varphi\L_\varphi$ enfin, nous aurons besoin d'un lemme classique qui peut être regardé comme la généralisation aux anneaux de polynômes sur $\Z_\varphi$ du Théorème 8 (iii) 
de \cite{se}~:
\end{itemize}

\begin{Lem} Soit $\Z_\varphi$ une extension finie non ramifiée de l'anneau $\Z_\l$ et $f_\varphi$ un polynôme distingué de l'anneau $\Z_\varphi[\g-1]$. Pour chaque entier naturel $n$ assez grand, il existe alors deux polynômes $a_n$ et $b_n$ dans $\Z_\varphi[\g-1]$ tels qu'on ait~:\smallskip

\centerline{$\sum_{k=0}^{\l-1} \g^{k\l^n} \ = \ \frac{\omega_{n+1}}{\omega_n} \ = \ 
\l (1+\l a_n) + b_n f_\varphi$ .}\smallskip

\noindent De plus l'élément $1+\l a_n$ est inversible dans l'anneau local $\L_\varphi= 
\Z_\varphi[[\gamma -1]]$.
\end{Lem}

\Preuve Raisonnons dans l'anneau quotient $\Z_\varphi[\gamma -1] / f_\varphi\Z_\varphi[\gamma -1]$. Si $d_\varphi$ est le degré du polynôme $f_\varphi$, nous avons par hypothèse~:

\centerline {$\ov{(\gamma -1)}^{\ d_\varphi} \equiv \ov0 \mod \ov \l, \qquad {\rm donc}\qquad \ov \gamma^{\,\l^{n-1}} -1 \equiv  \ov {(\gamma -1)}^{\,\l^{n-1}} \equiv \ov 1 \mod \ov \l$,}\smallskip

\noindent dès que $n$ est assez grand. Cela étant, il suit~:\smallskip

\centerline {$\ov \gamma ^{\,\l^{n-1}} \equiv \ov 1 \mod \ov \l ,\ {\rm puis} \ \
\ov \gamma ^{\,\l^n} \equiv \ov 1 \mod \ov \l^2 \ ;\ {\rm et \  enfin}\
\sum_{k=0}^{\l-1} \ov \g^{\,k\l^n} \equiv \ov \l \mod \ov \l^2 ,$}\smallskip

\noindent ce qui est précisément le résultat annoncé.\medskip

\noindent {\it Suite de la preuve du Théorème.} Appliqué au module isotypique $X = \L_\varphi / f_\varphi\L_\varphi$, le Lemme nous donne pour $n$ assez grand l'égalité $\nabla_{\!n+1}X = \l \,\nabla_{\!n} X$, donc~:\smallskip

\centerline{$(X:\nabla_{\!n}X) = (X:\nabla_{\!n_0}X) (\nabla_{\!n_0}X:\nabla_{\!n}X) = (X:\nabla_{\!n_0}X) (\nabla_{\!n_0}X:\l^{n-n_0}\nabla_{\!n_0}X)$}\smallskip

\noindent pour $n \gg n_0$ ; puis, comme attendu~:\smallskip

\centerline{$(X:\nabla_{\!n}X) \approx (\nabla_{\!n_0}X: \l^{n-n_0}\nabla_{\!n}X) = 
\l^{(n-n_0) d_\varphi \deg\varphi} \approx \l^n d_\varphi \deg \varphi$~;}\smallskip

\noindent ce qui achève la démonstration du cas élémentaire.\medskip

Dans la pratique, i.e. dans le contexte arithmétique de la théorie d'Iwasawa, il n'est pas possible en général d'appliquer directement le Théorème 1.4, car il est souvent nécessaire de faire appel à un résultat apparemment plus technique\footnote{C'est en particulier déjà le cas pour les groupes de classes d'idéaux au sens habituel.}~:

\begin{Cor}Soit $X$ un $\L[\D]$-module noethérien et de paramètres $\rho$, $\mu$ et $\lambda$~; et soit $Y$ un sous-module de $X$ de type fini sur $\Zl$. Alors, pour tout entier naturel $m$ assez grand, la suite des quotients définis pour $n \ge m$ par~:\smallskip

\centerline{$Y_n = X / (\nabla_{\!n} X + \frac{\omega_n}{\omega_m} Y)$}\smallskip

\noindent est paramétrée par les mêmes caractères $\rho$, $\mu$ et $\lambda$ que les $X_n = X/\nabla_{\!n} X$.
\end{Cor}

Enfin, il est évidemment  possible de jouer sur les deux termes qui interviennent dans la génération des idéaux $\nabla_{\!n}$. Il vient ainsi~:\medskip

\begin{Cor}Soit $X$ un $\L[\D]$-module noethérien de paramètres $\rho$, $\mu$ et $\lambda$. Si $\,\nabla_{\!n,k} = \l^{n+k} \L + \omega_n \L$ désigne l'idéal de l'algèbre d'Iwasawa $\L = \Z_\l[[\g-1]]$ engendré par l'élément $\l^{n+k}$ et le polynôme $\omega_n =\g^{\l^n}-1$, il existe un unique
caractère $\l$-adique virtuel $\nu$ du groupe $\D$ tel que l'ordre $\l^{x_n^\varphi}$ de la $\varphi$-composante du quotient $X/ \nabla_{\!n,k} X$ soit donné, pour chaque caractère $\l$-adique irréductible $\varphi$ et tout entier $n$ assez grand, par la formule~:\smallskip

\centerline{$x_n^\varphi \ = \ <\rho, \varphi>(n+k)\l^n \ + \ <\mu, \varphi>\l^n \ + \ <\lambda, \varphi>n \ + \ <\nu, \varphi>$.}\smallskip

\noindent En d'autres termes, la suite des quotients $\big(X / \nabla_{\!n,k} X\big)_{n 
\in \N}$ est paramétrée par~:\smallskip

\centerline{$\rho, \qquad \mu + (k-1) \rho , \quad et \quad \lambda$.}\smallskip
\end{Cor}

Toute la difficulté pour établir le Théorème des paramètres et ses corollaires est naturellement de vérifier que la formule annoncée est encore valide lorsque le module $X$ est seulement pseudo-isomorphe à un $\L[\D]$-module élémentaire.\smallskip

La démonstration de ce résultat fait l'objet de la section qui vient.

\bigskip

\noindent{\bf 1.3 Filtration associée aux idéaux $\nabla _{\!n}$ et preuve du Théorème}

\bigskip

Le Théorème des paramètres résulte d'une suite de quatre lemmes et du calcul déjà effectué de l'indice $(E : \nabla_{\!n} E)$ lorsque le module $E$ est élémentaire. 

\begin{Lem} Soit $E = T \oplus P$ un $\L[\D]$-module élémentaire somme directe de son sous-module de $\L$-torsion $T$ et d'un sous-module projectif $P$. Pour chaque entier naturel $n$, posons $\d^nE= \l^{n+1}E \cap \omega_n E$. Cela étant, pour tout $n \ge n_0$ assez grand, il vient~:
$$
\d^n E = \d^n T \oplus \d^n P \ ,\ avec\  \d^n T = \l^{n-n_0} \d^{n_0} T \  et \  
\d^n P = \frac{\omega_n}{\omega_{n_0}}\ \d^{n_0} P.
$$
\end{Lem}

\Preuve Observons d'abord que tout $\L[\D]$-module noethérien et projectif
$P$ étant facteur direct d'un module libre, nous avons immédiatement 
$\d ^nP = \l^{n+1} \cap \omega_nP$ par factorialité des composantes locales
$\L_\varphi$ de l'algèbre $\L[\D]$. En particulier, il suit comme annoncé~: 
$$
\d^nP = \frac{\omega_n}{\omega_{n_0}}\ \l^{n-n_0}\ \d^{n_0}P.
$$
Considérons maintenant le sous-module de torsion $T$. Ses facteurs 
indécompo\-sables, en nombre fini, sont de deux types~: des quotients de la
forme $\L_\varphi / \l^m \L_\varphi$, d'une part~; des quotients de la forme 
$\L_\varphi / f_\varphi \L_\varphi$, avec $f_\varphi$ distingué dans $\Z_\varphi[\g -1]$ 
d'autre part. Désignons par $T_0$ la somme directe des seconds. Pour $n$ 
assez grand, nous avons $\l^n T = \l^n T_0$, donc $\d^n T = \d^n T_0$. 
Appliquons maintenant le Lemme 2.2 à chacun des facteurs 
indécomposables de $T_0$. Pour $n$ assez grand, disons pour $n \ge n_0$,
nous obtenons $\frac{\omega_{n+1}}{\omega_n} T_0 = \l T_0$, d'où~:
$$
\d^n T_0 = \l^{n+1} T_0 \cap \omega_n T_0 = \l^{n-n_0}\l^{n_0+1} T_0 \cap 
\frac{\omega_n}{\omega_{n_0}}\ \omega_{n_0} T_0 = \l^{n-n_0}
(\l^{n_0+1} T_0 \cap \omega_{n_0} T_0).
$$
Finalement il vient bien, comme attendu~: $\d^n T_0 = \l^{n-n_0} \d^{n_0} T_0$.

\begin{Lem} Soit $X$ un sous-module d'indice fini d'un module
élémentaire $E$. Pour chaque caractère $\l$-adique irréductible $\varphi$, 
la suite $\big((\d^n Y_\varphi : \d^n X_\varphi)\big)_{n \in \N}$ des $\varphi$-parties des 
indices $(\d^n E_\varphi : \d^n X_\varphi)$ est stationnaire.
\end{Lem}

\Preuve Désignons par $E^{tor}$ le sous-module de $\L$-torsion de $E$  et notons
$p$ la projection canonique $E \to E/E^{tor}$. Nous avons alors~:
$$
(\d^n E_\varphi : \d^n X_\varphi) = \frac{(\d^n E_\varphi : \d^n E^{tor}_\varphi + \d^n X_\varphi)}
{(\d^n E^{tor}_\varphi + \d^n X_\varphi : \d^n X_\varphi)} = \frac{(p(\d^n E_\varphi) :
 p(\d^n X_\varphi))}{(\d^n E^{tor}_\varphi : \d^n X^{tor}_\varphi)} \ .
$$
Et il s'agit de voir que numérateur et dénominateur sont tous deux stationnaires.
Or d'un côté, nous avons  l'égalité évidente $p(\d^{n+1}E_\varphi) = \l \ 
\frac{\omega_{n+1}}{\omega_n} \ p(\d^{n}E_\varphi)$ ainsi que l'inclusion immédiate
$p(\d^{n+1}X_\varphi) \subset \l \ \frac{\omega_{n+1}}{\omega_n} \ p(\d^{n}X_\varphi)$~;
d'où il résulte que le numérateur, qui va donc décroissant pour $n$ assez 
grand, est bien stationnaire.\smallskip

Et d'un autre côté, nous avons $\d^{n+1} E^{tor}_\varphi =\l \ \d^n E^{tor}_\varphi$ et
 $\d^{n+1} X^{tor}_\varphi =\l \ \d^n X^{tor}_\varphi$ pour $n \ge n_0$, par une extension
facile du Lemme 2.3, puisque $E^{tor}_\varphi$ et $X^{tor}_\varphi$  sont annulés
par un même polynôme distingué $f_\varphi \in \Z_\varphi[\g-1]$. Cela étant, comme
pour $n_0$ assez grand la multiplication par $\l^{n_0}$ a tué le sous-module 
de $\Zl$-torsion de $E^{tor}_\varphi$, la multiplication par $\l$ est injective dans 
$\l^{n_0}E^{tor}_\varphi$ et le dénominateur est ainsi stationnaire pour 
$n \ge n_0$.

\begin{Lem}Soit $X$ un sous-module d'indice fini d'un $\L[\D]$-module élémentaire $E$. Pour chaque caractère $\l$-adique irréductible $\varphi$, la suite $\big((\nabla_{\!n} Y_\varphi : \nabla_{\!n} X_\varphi)\big)_{n \in \N}$ des $\varphi$-parties des indices $(\nabla_{\!n} E_\varphi : \nabla_{\!n} X_\varphi)$ est stationnaire.
\end{Lem}

\Preuve Dans la suite exacte courte de $\Zl[\D]$-modules finis
$$
0 \rightarrow \d^n E_\varphi / \d^n X_\varphi \rightarrow ( \l^n E_\varphi / \l^n X_\varphi) 
\oplus (\omega_n E_\varphi / \omega_n X_\varphi) \rightarrow
\nabla_{\!n} E_\varphi / \nabla_{\!n} X_\varphi \rightarrow 0
$$
les termes médians, qui vont décroissant, sont stationnaires. Le résultat
annoncé découle donc directement du Lemme 1.10 ci-dessus.

\begin{Lem} Soit $X$ un $\L[\D]$-module noethérien et $E$ l'unique 
$\L[\D]$-module élémentaire auquel il est pseudo-isomorphe. Pour chaque
caractère $\l$-adique irréductible $\varphi$, la suite des quotients $\big( X_\varphi 
: \nabla_{\!n} X_\varphi \big) / \big( E_\varphi : \nabla_{\!n} E_\varphi \big)$ est stationnaire.
\end{Lem}

\Preuve Prenons un pseudo-isomorphisme $t$ de $X$ vers $E$ et considérons
la suite exacte courte associée~:
$$\CD
0 @>>> N @>>> X @>t>> E @>>> C @>>> 0.
\endCD$$
Le noyau $N$ et le conoyau $C$ étant finis, choisissons $n_0$ assez grand
pour avoir $N \cap \nabla_{\!n_0}X = 0$ et $\nabla_{\!n_0}C = 0$, i.e. 
$\nabla_{\!n_0}E \subset t(X)$. Pour chaque $n \ge n_0$, nous obtenons alors
la suite exacte~:
$$
0 \rightarrow {}^{-1}t(\nabla_{\!n}E)/ \nabla_{\!n}X \rightarrow X/ \nabla_{\!n}X 
\rightarrow E/ \nabla_{\!n}E \rightarrow C \rightarrow 0\ ;
$$
où, pour chaque caractère irréductible $\varphi$ du groupe $\Delta$, l'ordre de la 
$\varphi$-partie du terme de gauche est donnée par la formule~:
$$
({}^{-1}t(\nabla_{\!n}E_\varphi) : \nabla_{\!n}X_\varphi) = | N_\varphi | \  (\nabla_{\!n}E_\varphi) : 
\nabla_{\!n}t(X_\varphi)).
$$
Elle est donc ultimement constante  en vertu du Lemme 1.11.
\medskip

Nous pouvons maintenant achever la démonstration des résultats annoncés dans la section 2.\medskip

\noindent {\it Preuve du Théorème des paramètres et de ses corollaires.} Distinguons les cas~:\smallskip 

{\it (i) } Le Théorème 1.3 résulte directement  du calcul déjà effectué de l'indice $\big( E_\varphi : \nabla_{\!n} E_\varphi \big)$ dans le cas élémentaire et du Lemme 1.12 ci-dessus.\smallskip

{\it (ii)} Le Corollaire 1.7 s'en déduit aisément~: compte tenu de la suite exacte courte canonique qui relie les $X_n$ avec les $Y_n$~:
$$
1 \rightarrow  \frac{\omega_n}{\omega_m} Y / (\frac{\omega_n}{\omega_m} Y \cap \nabla_{\!n} X) \rightarrow X_n = X / \nabla_{\!n} X \rightarrow Y_n = X / (\nabla_{\!n} X + \frac{\omega_n}{\omega_m} Y) \rightarrow 1,
$$
tout revient à vérifier que les noyaux à gauche sont ultimement constants. Or cela est clair, puisqu'ils vont naturellement décroissant et qu'ils sont évidemment finis car de type fini sur $\Zl$ par hypothèse et d'exposant fini par construction.\smallskip

{\it (iii)} Enfin, le Corollaire 1.8 s'obtient de la même façon que le Théorème en remplaçant les idéaux $\nabla_{\!n}$ par les idéaux $\nabla_{\!n, k}$ pour un $k$ fixé~: les calculs étant identiques, il est possible d'écrire la même série de lemmes que ci-dessus, ce qui ramène la preuve du Corollaire au cas élémentaire, lequel est immédiat.

\newpage


\setcounter{section}{2}  \setcounter{Th}{0} 

\noindent{\large \bf 2. Applications à l'arithmétique des tours cyclotomiques}


\bigskip

Nous supposons désormais que le nombre premier $\l$ est impair.\medskip

Dans cette section, nous désignons par $F$ un corps de nombres totalement réel et par $K$ une extension abélienne de $F$, de groupe de Galois $\Delta\, =\, \Gal (L/F)$, de degré $d$ étranger à $\ell$, contenant une racine primitive $\ell$-ième de l'unité $\zeta$.\smallskip 

L'hypothèse d'imparité $\ell \ne 2$ nous assure que le corps $K$ est totalement imaginaire et que le groupe de Galois $\D$ contient la conjugaison complexe $\tau$. Dans ce contexte, il est naturel d'appeler {\it réels} les caractères irréductibles $\varphi$ de $\D$ qui vérifient l'inégalité $\varphi(\tau) > 0$~;  et {\it imaginaires} ceux qui vérifient l'identité opposée $\varphi(\tau) <0$. Plus généralement~:

\begin{Def} {\rm Nous disons qu'un caractère $\l$-adique virtuel du groupe $\D$ est~:}\par

(i) totalement réel, {\rm lorsque tous ses facteurs irréductibles sont réels, i.e. lorsque tous ses facteurs absolument irréductibles prennent la valeur +1 sur la conjugaison complexe $\tau$~;}\par

(ii) totalement imaginaire, {\rm lorsque tous ses facteurs irréductibles sont imaginaires, i.e. lorsque tous ses facteurs absolument irréductibles prennent la valeur -1 sur la conjugaison complexe $\tau$.}
\end{Def}

Sommant sur toutes les composantes isotypiques réelles de l'algèbre $\Zl[\D]$, on définit en particulier sa composante réelle, qui n'est rien d'autre que l'image de $\Zl[\D]$ par l'idempotent $e_\oplus = \frac{1}{2}(1+\tau)$~:\smallskip

\centerline{$\Zl[\D]^\oplus =  \Zl[\D] e_\oplus = \oplus_{\varphi(\tau)>0}\Zl[\D] e_\varphi $~;}\smallskip

\noindent et on définit également sa composante imaginaire à partir cette fois de l'idempotent complémentaire $e_\ominus = \frac{1}{2}(1-\tau)$ par~:\smallskip

\centerline{$\Zl[\D]^\ominus =  \Zl[\D] e_\ominus = \oplus_{\varphi(\tau)<0}\Zl[\D] e_\varphi$~;}
\smallskip

\noindent les mêmes décompositions valant, plus généralement, pour tout $\Zl[\D]$-module $M$. En particulier, il est commode de décomposer chaque caractère $\l$-adique virtuel $\chi$ du groupe $\D$ sous la forme~:\smallskip

\centerline{$\chi = \chi^\oplus + \chi^\ominus$,}\smallskip

\noindent en regroupant séparément facteurs réels et facteurs imaginaires.\medskip

Parmi les caractères de $\Delta $ figurent en particulier le {\it caractére unité} 1, dont l'idempotent associé est donné à partir de la norme algébrique $\nu _\Delta \, =\,\sum_{\delta\in \Delta}\delta $ par 
$e_1\, =\, \frac{1}{d}\ \nu _{\Delta }$, et le {\it caractère cyclotomique} $\omega$, caractérisé par 
l'identité~:\smallskip

\centerline{$\zeta ^\sigma\, =\, \zeta ^{\omega (\sigma )}\quad \forall\, \sigma\, \in \Delta$.}\smallskip

\noindent C'est aussi le caractère du module de Tate $\Tl = \varprojlim \mu(K_n)$ construit sur les $\l$-groupes de racines de l'unité.\smallskip

L'inverse\footnote{ De façon générale, il est commode de noter $\psi ^{-1}$ le caractère $\sigma\mapsto \psi (\sigma ^{-1})$} $\bar\omega =\omega ^{-1}$ de $\omega $, c'est à 
dire le caractère défini par $\bar\omega (\sigma )=\omega (\sigma ^{-1})$, est dit souvent $anticyclotomique$, et l'involution\smallskip

\centerline{$\psi\mapsto\psi ^* =\omega\psi ^{-1}$,}\smallskip

\noindent de l'algèbre $R_{\Z_{\ell} }(\Delta )$ des caractères $\ell$-adiques virtuels de $\Delta $ est connue traditionnellement sous le nom d'{\it involution du miroir}~; on dit encore que $\psi ^* $ est le {\it reflet} de $\psi$. En particulier, le caractère cyclotomique $\omega = 1^*$ étant imaginaire, il suit que le reflet d'un caractère réel est imaginaire et vice versa.

\newpage

\noindent{\bf 2.1 Le contexte arithmétique de la Théorie d'Iwasawa}

\medskip

Introduisons maintenant la $\Zl$-extension cyclotomique $F_\infty = \bigcup_{n\in \N} F_n$ du corps totalement réel $F$, en convenant de désigner par $F_n$ l'unique sous-extension de $F_\infty$ qui est de degré $\l^n$ sur $F$ (de sorte que l'on a $F_0 = F$)~; notons $\Gamma = \gamma^{\Z_\l}$ le groupe de Galois $\Gal(F_\infty
/F)$ identifié à $\Zl$ par le choix d'un générateur topologique $\gamma$~; et écrivons $\Lambda = \Z_\l[[\gamma - 1]]$  l'algèbre d'Iwasawa associée.\par

Notons de même $K_\infty = \bigcup_{n\in \N} K_n$ la $\Zl$-extension cyclotomique de $K$ en désignant par $K_n = KF_n$ le compositum de $K$ et de $F_n$. L'extension $K_\infty/F$ est encore abélienne~; son groupe de Galois $\Gal (K_\infty/F)$ est le produit $\Gamma \times \D$~;  et l'algèbre complète associée est l'algèbre $\Zl[[\gamma-1]][\D] =\L[\D]$ étudiée plus haut.\smallskip

Fixons alors deux ensembles finis disjoints $S$ et $T$ de places de $K$ stables par $\D$~; notons $L$ l'ensemble des places $\l$-adiques (i.e. au-dessus de $\l$) de $K$~; puis $S^\l = S \cap L$ (respectivement $T^\l  = T \cap L$) la partie sauvage de $S$ (resp. $T$) et $S^0 = S \setminus S^\l$ 
(resp. $T^0  = T \setminus T^\l$) la partie modérée. Enfin, pour chaque indice $n \in \N \cup \{\infty\}$, notons de même $S_n = S_n^\l \cup S_n^0$ (resp. $T_n = T_n^\l \cup T_n^0$) l'ensemble (fini) des places de $K_n$ qui sont au-dessus de $S$ (resp. $T$).\smallskip

Nous sommes intéressés dans ce qui suit par le groupe de Galois $\C_{S_\infty}^{T_\infty}$ de la pro-$\l$-extension abélienne $T_\infty$-ramifiée $S_\infty$-décomposée (i.e. non ramifiée en dehors des places au dessus de $T$ et complètement décomposée aux places au dessus de $S$) maximale $H_{S_\infty}^{T_\infty}$ de $K_\infty$~; autrement dit par la limite projective $\varprojlim \C_{S_n}^{T_n}$ des mêmes groupes attachés aux étages finis $K_n$ de la tour cyclotomique, groupes que la théorie $\l$-adique du corps de classes identifie aux $\l$-groupes de $S$-classes $T$-infinitésimales des corps $K_n$, c'est à dire aux quotients~:\smallskip

\centerline{$ \Cl_{S_n}^{T_n} = \J_{K_n} / \J_{S_n}^{T_n} \R_{K_n} $,}\smallskip

\noindent où $\J_{K_n} = \prod^{\rm res}\R_{\p_n} $ désigne le $\l$-adifié du groupe des idèles du corps $K_n$~; $ \J_{S_n}^{T_n}=\prod_{\p_n \notin T_n}\U_{\p_n} \prod_{\p_n \in S_n}\R_{\p_n}$ le sous-groupe des $S$-idèles $T$-infinitésimaux~; et $\R_{K_n}= \Zl \otimes K_n^\times$ le sous-groupe principal de $\J_{K_n}$ (cf. \cite{ja2}).\par

Avec les conventions ci-dessus et celles de la section 1, nous avons~:

\begin{Th} Le groupe $\C_{S_\infty}^{T_\infty}$ est un $\L[\D]$-module noethérien isomorphe à la limite projective $\varprojlim \Cl_{S_n}^{T_n}$ des $\l$-groupes de $S$-classes $T$-infinitésimales des corps $K_n$ ainsi qu'à celle de leurs quotients respectifs d'exposant $\l^{n+1}$~:\smallskip

\centerline{$\C_{S_\infty}^{T_\infty}  \simeq \varprojlim \,{}^{\l^{n+1}}\! \Cl_{S_n}^{T_n}$ .}\smallskip 

\noindent Et les trois caractères $\rho^T_S$, $\mu^T_S$ et $\lambda^T_S$, qui lui sont associés paramètrent aussi la suite des quotients finis ${}^{\l^{n+1}}\! \Cl_{S_n}^{T_n}$, à l'exception du cas spécial $S=\emptyset$ et $T=L$ dans lequel le paramètre lambda doit être augmenté du caractère unité.
\end{Th}

\Preuve Commençons par établir proprement ce premier résultat. Ecrivons pour simplifier $X$ pour $\,\C_{S_\infty}^{T_\infty}$ et $X_n$ pour $\,{}^{\l^{n+1}}\!\Cl_{S_n}^{T_n}$. Observons d'abord que les groupes $X_n$ ne sont pas {\it tout à fait} les quotients $X/ \nabla_{\!n}X$, mais qu'ils en diffèrent d'une quantité bornée~: en effet, le groupe de Galois $\Gamma_n =\Gal (K_\infty/K_n) = \gamma^{\l^n\Zl}$ étant procyclique, le plus grand quotient de $X$ fixé par $\Gamma_n$, i.e. le groupe $ {}^{\Gamma_n}\!X = X/ \omega_n X$, correspond à la plus grande sous-extension $H_S^T(K_\infty /K_n)/K_\infty$ de l'extension $H_S^T(K_\infty)/K_\infty$ qui provient par composition avec $K_\infty$ d'une pro-extension abélienne de $K_n$. En particulier, une telle extension étant non ramifiée en dehors de $L$ et de $T$, les quotients $X/ \nabla_{\!n}X$ sont finis, en vertu du corps de classes, ce qui montre, d'après la proposition 1.3, que $X$ est bien un $\L[\D]$-module noethérien.\par

Plus précisément, fixons $m$ assez grand pour que les places au-dessus de $S$ soient totalement inertes et celles au-dessus de $L\setminus T^\l$ totalement ramifiées dans l'extension $K_\infty/K_n$~; et notons $Y_m$ le sous-module de $X$, de type fini sur $\Zl$, qui est engendré par les sous-groupes (procycliques) de décomposition des premières et les sous-groupes d'inertie  des secondes. Il vient alors~:\smallskip

\centerline{${}^{\l^{m+1}}\!\Gal(H_S^T(K_\infty/K_m)/K_\infty) = X / (\nabla_{\!n} X +  Y_m)$~;}\smallskip

\noindent et par un argument classique de théorie d'Iwasawa (cf. e.g. \cite{la}, \cite{se}, ou \cite{wa} dans le cas classique des groupes de classes ordinaires)~:\smallskip

\centerline{${}^{\l^{n+1}}\!\Gal(H_S^T(K_\infty/K_n)/K_\infty) = X / (\nabla_{\!n} X + \frac{\omega_n}{\omega_m}  Y_m)$.}\smallskip

Le résultat annoncé en résulte en vertu du Corollaire 1.7. puisqu'en dehors du cas spécial on a directement par l'isomorphisme du corps de classes~:\smallskip

\centerline{${}^{\l^{n+1}}\!\Cl_S^T(K_n) \simeq {}^{\l^{n+1}}\!\Gal(H_S^T(K_\infty/K_n)/K_\infty)$~;}\smallskip

\noindent tandis que dans le cas spécial, il vient~:\smallskip

\centerline{${}^{\l^{n+1}}\!\Cl_S^T(K_n) \simeq {}^{\l^{n+1}}\!\Gal(H_S^T(K_\infty/K_n)/K_\infty) \oplus \Gamma/ \Gamma_n$~;}\smallskip

\noindent puisqu'il convient alors de tenir compte en outre de la contribution de la tour cyclotomique $K_\infty/K_n$, lorsque celle-ci est elle-même $S$-décomposé et $T$-ramifiée, c'est à dire lorsqu'on a $S=\emptyset$ et $T=L$.\medskip

Nous nous proposons dans ce qui suit d'étudier les trois caractères définis ci-dessus à la lumière des égalités du miroir obtenues par  Georges Gras (cf. \cite{gr1}). 
Dans notre contexte, ces identités relient les invariants d'Iwasawa relatifs à un couple $(S,T)$ aux mêmes invariants relatifs au couple transposé $(T,S)$. Elles s'obtiennent en comparant les informations données par la théorie du corps de classes d'une part, et la théorie de Kummer d'autre part en présence des racines  de l'unité. C'est pourquoi , nous ne considérons dans cette section que le cas de la $\Zl$-extension {\it cyclotomique} $K_\infty = K^c$ d'un corps $K$ {\it supposé déjà contenir} les racines $\l$-ièmes de l'unité~: en l'absence de l'une ou l'autre de ces deux hypothèses, nos arguments de dualité ne pourraient s'appliquer.\par

Cela dit, fixons $n$ assez grand pour que les places de $S \cup T$ ne se décomposent pas dans l'extension $K_\infty /K_n$ et plaçons nous au niveau $n$ de la tour en omettant occasionnellement l'indice $K_n$ dans ce qui suit. Les résultats de \cite{jm} se généralisent alors comme suit~:

\begin{Prop} Avec les conventions précédentes, lorsque la réunion $S \cup T$ contient l'ensemble $L$ des places $\l$-adiques, le radical kummérien $\Rad(H^{T_n}_{S_n}/K_n)$ attaché à la $\l$-extension abélienne maximale $H^{T_n}_{S_n}$ de $K_n$ qui est $S$-décomposée, $T$-ramifiée et d'exposant $\l^{n+1}$ est donné par~:\smallskip

$ $(i) \ $\Rad^{S_n}_{T_n} = \{ \l^{-(n+1)} \otimes x \in \l^{-(n+1)} \Zl / \Zl \otimes \R_{K_n}\, | \ x \in \J^{S_n}_{T_n, K_n} \J_{K_n}^{\l^{n+1}}\}$~;\smallskip

\noindent et le groupe de Galois $\Gal(H^{T_n}_{S_n}/K_n)$ de la même extension $H^{T_n}_{S_n}$ est donné par~:\smallskip

(ii) \ $\Gal^{T_n}_{S_n} \simeq\ {}^{\l^{n+1}}\!\Cl^{T_n}_{S_n}(K_n) =\ {}^{\l^{n+1}}\!( \J_{K_n} / \J^{T_n}_{S_n} \R_{K_n}) = \J_{K_n} / \J^{T_n}_{S_n}  \R_{K_n} \J_{K_n}^{\l^{n+1}}$.\smallskip

\noindent Ces deux modules se correspondant dans la dualité de Kummer, on a en outre~:\smallskip

\centerline{$\Gal^{T_n}_{S_n} \simeq \Hom _{\Zl[\D]} (\Rad^{S_n}_{T_n}, {}_{\l^{n+1}} \mu_{\l_\infty})\ $ et $\  \Gal^{T_n}_{S_n} \simeq \Hom _{\Zl[\D]}  (\Rad^{S_n}_{T_n},{}_{\l^{n+1}}\mu_{\l^\infty})$.}\smallskip

En particulier, en dehors du cas spécial $S=\emptyset$ et $T=L$ dans lequel le paramètre lambda doit être augmenté du caractère unité, la suite $(\Gal^{T_n}_{S_n})_{n\in\N}$ est paramétrée par le triplet de caractères $(\rho^T_S, \mu^T_S, \lambda^T_S)$~; et la suite $(\Rad^{T_n}_{S_n})_{n\in\N}$  par le triplet $(\rho^{T*}_S, \mu^{T*}_S, \lambda^{T*}_S)$, image du précédent dans l'involution du miroir. 
\end{Prop}

\Preuve Il suffit d'exprimer, en termes kummériens dans le premier cas, de corps de classes dans le second, les conditions locales qui définissent l'extension.\medskip

\bigskip

\noindent{\bf 2.2 Isomorphismes du miroir dans les tours cyclotomiques}

\bigskip

La clef du Spiegelungssatz de Gras dans notre situation consiste à introduire un groupe de congruences {\it ad hoc}, qui n'est ni un radical ni un groupe de classes, mais dont l'écart avec l'un et l'autre est parfaitement maîtrisé.\par

Supposons donc $T$ non vide puis, pour chaque place $\p_n$ de $T_n$, choisissons un entier 
$\nu_{\p_n}$ assez grand pour que le groupe supérieur d'unités locales $\U_{\p_n^{\nu_{\p_n}}}$ soit sans torsion et contenu dans la puissance $\l^{n+1}$-ième du groupe principal $\U_{\p_n}$~; cela fait, introduisons le diviseur $\m_n \ =\ \prod_{\p_n \in T_n} \p_n^{\nu_{\p_n}}$~; notons $\R^{\m_n}$ le groupe des idèles principaux de $K_n$ congrus à 1 modulo $\m_n$, et $\E_{S_n}^{\m_n}$ le sous-groupe de $\R^{\m_n}$ construit sur les $S$-unités. Définissons alors le {\it pseudo-radical} ${}_{\l^{n+1}} \PR^{\m_n}_{S_n}$  par~:\smallskip

\centerline{${}_{\l^{n+1}} \PR^{\m_n}_{S_n}   = \{\l^{-(n+1)}\! \otimes x \in \l^{-(n+1)} \Zl / \Zl \otimes \R^{\m_n} | \ x \in \J_{S_n} \J^{\l^{n+1}}\}$,}\medskip

\noindent où $\J_{S_n} = \prod_{\p_n \in S_n} \R_{\p_n} \prod_{\p_n \notin S_n} \U_{\p_n}$ est le groupe des $S$-idèles~; notons de même~: \smallskip

\centerline{${}_{\!\l^{n+1}} \EE^{\m_n}_{S_n} = \{\l^{-(n+1)}\! \otimes x \in \l^{-(n+1)} \Zl / \Zl \otimes \R^{\m_n} | \ x \in \E_{S_n} ^{\m_n} \} \simeq {}^{\l^{n+1}}\E^{\m_n}_{S_n}$}\smallskip

\noindent son sous-groupe construit sur les $S$-unités de $\E_{S_n}^{\m_n}$~; et écrivons enfin ${}_{\l^{n+1}}\mu$ le $\l$-groupe des racines $\l^{n+1}$-ièmes de l'unité. Nous obtenons comme dans \cite{jm}~:

\begin{Lem} Avec ces notations, le pseudo-radical ${}_{\l^{n+1}} \PR^{\m_n}_{S_n}$ est relié aux groupes précédents par les suites exactes courtes de modules galoisiens~:
$$ \begin{CD}
1@>>>{}_{\l^{n+1}}( \R / \R^{\m_n})/{}_{\l^{n+1}}\mu@>>>{}_{\l^{n+1}}\PR^{\m_n}_{S_n}
@>\phi_n>>{}_{\l^{n+1}}\!\Rad^{T_n}_{S_n}@>>>1\\ 
@. @.  @\vert  @.\\ 
1 @>>> {}_{\l^{n+1}} {\EE}^{\m_n}_{S_n} @>>> {}_{\l^{n+1}}\PR^{\m_n}_{S_n}@>\psi_n>> 
{}^{\l^{n+1}}\!\Gal^{T_n}_{S_n}  @>>>1\ ,
\end{CD} $$
\end{Lem}

La première suite exacte n'est autre que la suite exacte canonique donnée dans \cite{jm}. La seconde appelle une petite explication~: {\it a priori}, le morphisme canonique $\psi_n$ prend ses valeurs dans le sous-groupe de $\l^{n+1}\!$-torsion du $\l$-groupe $Cl_{S_n}^{\m_n}$ des $S$-classes de rayon $\m_n$~; mais, celui-ci étant fini, nous avons bien~:\smallskip

\centerline{${}_{\l^{n+1}}\Cl_{S_n}^{\m_n} \ \simeq\ {}^{\l^{n+1}}\!\Cl_{S_n}^{\m_n} \ =\ {}^{\l^{n+1}}\!\Cl_{S_n}^{T_n}\ \simeq \ {}^{\l^{n+1}}\!\Gal_{S_n}^{T_n}$.}\medskip

Dans le diagramme obtenu, les deux modules à droite sont paramétrés et leurs caractères invariants sont donnés par la Proposition 2.3 ci-dessus. Pour calculer ceux des modules de gauche, introduisons quelques notations~:\smallskip

Pour chaque place $\p_n$ du corps $F_n$, notons $\Delta_\p$ le sous-groupe de décomposition de $\p_n$ dans l'extension abélienne $K_n/F_n$ (qui ne dépend que de la place $\p$ de $F$ au-dessous de $\p_n$)~; notons $\chi_\p$ l'induit à $\Delta$ du caractère unité de $\Delta_\p$~; et pour tout ensemble fini $T$ de places finies définissons le caractère $\chi_T$ par~:

\centerline{$\chi_T = \sum_{\p_\infty \in T_\infty} \chi_\p$,}\medskip

\noindent où la somme porte sur les places de $F_\infty$ au-dessus de $T$ (lesquelles sont bien en nombre fini). Définissons alors le degré $\l$-adique total en $T$ par la formule~:\smallskip

\centerline{$\deg_\l T = \sum_{\p \in T_\l} [F_\p : \Ql]$,}\smallskip

\noindent où la somme porte cette fois sur les places $\l$-adiques de $F$ contenues dans $T$. Et posons enfin (en notant $\chi_{\rm reg}$ le caractère régulier du groupe $\D$)~:\smallskip

\centerline{$\delta_T = \deg_\l T\ \chi_{\rm reg}$.}\smallskip

\noindent Posons enfin~: $\chi_\infty =[F:\Q]\ \chi_{\rm reg}^\oplus$. Cela étant, nous avons~:

\begin{Lem}  Notons $\U_{T_n} = \prod_{\p_n \in T_n} \U_{\p_n}$ la $T$-partie du goupe des idèles unité de $K_n$ et $\U_{\m_n}$ le sous-groupe de congruence attaché à $\m_n$. \smallskip

(i) L'hypothèse $\U_{\m_n} \subset \U_{T_n}^{\l^{n+1}}$ et le théorème d'approximation faible donnant~:\smallskip

\centerline{${}_{\l^{n+1}}( \R / \R^{\m_n}) \simeq {}_{\l^{n+1}}(\U_{T_n} / \U_{\m_n})  
\simeq {}^{\l^{n+1}}\!(\U_{T_n} / \U_{\m_n}) =  {}^{\l^{n+1}}\!\U_{T_n}$~;}\smallskip

\noindent les modules ${}_{\l^{n+1}}( \R / \R^{\m_n})/{}_{\l^{n+1}}\mu$ sont paramétrés par le triplet $(\delta_T, 0, (\chi_T-1)^*)$.\smallskip

(ii) De même, si $\E_{S_n}^{tor}$ est le $\l$-groupe des racines de l'unité dans $K_n$, il vient~:\smallskip

\centerline{${}_{\l^{n+1}} \EE^{\m_n}_{S_n}\ \simeq\ {}^{\l^{n+1}}\E^{\m_n}_{S_n}\ 
\simeq\ {}^{\l^{n+1}}\!(\E_{S_n}/ \E_{S_n}^{tor})$~;}\smallskip

\noindent et les modules ${}_{\l^{n+1}} \EE^{\m_n}_{S_n}$ sont ainsi paramétrés par le triplet $(\chi_\infty, 0, \chi_S-1)$.
\end{Lem}

\Preuve Il s'agit de reécrire la proposition 5 de \cite{jm} en termes de caractères~:\smallskip

Pour $(i)$, il suffit d'observer d'une part que le quotient sans torsion de $\U_{T_n}$ est un $\Zl[\D]$-module projectif de carctère $[F_n :\Q]\ \chi_{\rm reg}$, ce qui donne l'expression du paramètre $\rho$~; et d'observer d'autre part que son sous-module de $\l^{n+1}$-torsion est un $\Z / \l^{n+1}\Z[\D]$-module de caractère $\omega \chi_T = \chi_T^*$, ce qui donne l'invariant $\lambda$.\smallskip

Pour $(ii)$, cela revient à remplacer le théorème de Dirichlet sur le rang des $S$-unités par le théorème de représentation de Herbrand pour ces mêmes groupes.\medskip

Rassemblant ces résultats, nous obtenons alors le Théorème du miroir~:

\begin{Th}Si $S \cup T$ contient l'ensemble $L$ des places $\l$-adiques, les paramètres d'Iwasawa des $\l$-groupes $^{\l^{n+1}}\!Cl ^T_S(K_n)$ vérifient les identités du miroir~:
\begin{enumerate}
\item[(i)] $\rho_S^T +\frac{1}{2} (\chi_\infty + \delta_S) =
[\rho_T^S +\frac{1}{2} (\chi _\infty + \delta_T)]^*$ ;
\item[(ii)] $\mu_S^T=\mu_T^{S*}$ ;
\item[(iii)] $\lambda_S^T + (\chi _S - 1 ) = [\lambda_T^S + (\chi_T -1)]^*$.
\end{enumerate}
\end{Th}

\Preuve Les identités $(ii)$ et $(iii)$ s'obtiennent immédiatement en exprimant de deux façons les paramètres des pseudo-radicaux ${}_{\l^{n+1}}\PR^{\m_n}_{S_n}$ à l'aide de la Proposition 2.3 et du Lemme 2.5. Dans le cas de $(i)$, le même calcul donne directement~:\smallskip

\centerline{$\rho^T_S + \chi_\infty = \rho_T^{S*} + \delta_T$,}\smallskip

\noindent d'où l'on déduit la formule annoncée en écrivant l'identité reflet et en observant que les deux caractères $\delta_T$ et $\delta_S$ sont invariants par le miroir et de somme~:\smallskip

\centerline{$\delta_T +\delta_S = [F : \Q]\ \chi_{\rm reg} = \chi_\infty + \chi_\infty^*$.}\medskip

Tirons une première conséquence immédiate de ce résultat de dualité~:

\begin{Th} Le paramètre $\rho^T_S$ est donné pour $S$ et $T$ arbitraires par~:\smallskip

\centerline{$\rho^T_S = \delta_T^\ominus$.}
\end {Th}

\Preuve Observons d'abord que le paramètre $\rho_S^T$ est bien indépendant de l'ensemble $S$, puisque le sous-module de $\C^T(K_\infty)$ qui est engendré par les sous-groupes de décomposition respectifs des places de $S$ est de rang fini sur $\Zl$. Quitte à grossir $S$, nous pouvons donc supposer, sans restreindre la généralité, que $S \cup T$ contient l'ensemble $L$ de toutes les places $\l$-adiques. \par

Le résultat annoncé résulte alors tout simplement du fait que le paramètre $\rho^T_S$ est totalement imaginaire, lequel est lui même une conséquence du fait que le défaut de la conjecture de Leopoldt dans la tour cyclotomique $K_\infty^+/K^+$ construite sur le sous-corps totalement réel $K^+$ de $K$ (i.e. le nombre de $\Zl$-extensions linéairement indépendantes de ses étages finis $K_n^+$) est borné. Il vient donc bien~:\smallskip

\centerline{$ \rho^T_S = \rho^{T\ominus}_S = (\rho^T_S + \chi_\infty)^\ominus = (\rho_T^{S *} + \delta_T)^\ominus = (\rho_T^{S\oplus})^* + \delta_T^\ominus =\delta_T^\ominus$}

\bigskip

\noindent{\bf 2.3 Application au calcul des paramètres $\mu_S^T$ et $\lambda_S^T$}

\bigskip

Intéressons nous d'abord aux invariants $\mu$. Observons en premier lieu que le paramètre $\mu_S^T$ est en fait indépendant de l'ensemble $S$ (puisque le sous-module de $\C^T(K_\infty)$ qui est engendré par les sous-groupes de décomposition respectifs des places de $S$ est de rang fini sur $\Zl$) comme de la partie modérée $T^0$de l'ensemble $T$ (puisqu'il en est de même du sous-module engendré par les sous-groupes de décomposition des places de $T^0$). En d'autres termes, nous avons~:\smallskip

\centerline{$\mu^T_S = \mu^{T^\l}$,}\smallskip

\noindent pour tous $S$ et $T$. Compte tenu des identités du miroir, il vient ainsi~:

\begin{Th} Soient $S$ et $T$ quelconques~; puis $T^\l$ la partie sauvage de $T$ (i.e. le sous-ensemble des places $\l$-adiques de $T$) et $\ov{T}^\l = L \setminus T^\l$ son complémentaire dans L. Avec ces notations, le paramètre $\mu^T_S$ est donné par la formule~:\smallskip

\centerline{$\mu^T_S = \mu^{T^\l} = \mu^{T^\l}_{\ov {T}^\l} = \mu_{T^\l}^{\ov {T}^\l  *} = \mu^{\ov {T}^\l *}$.}\smallskip
\end{Th}

Et l'on voit que la détermination des divers invariants $\mu^T_S$ se ramène ainsi à celle des seuls paramètres $\mu^{T^\l}$, où $T^\l$ décrit les sous-ensembles de $L$ stables par $\D$. Plus précisément, il suffit même de ne considérer que les invariants réels~:\smallskip

\begin{Cor} Avec les notations du Théorème, on a par dualité~:\smallskip

\centerline{$\mu^T_S = \mu^{T^\l} = \mu^{T^\l \oplus} +  (\mu^{\ov {T}^\l \oplus})^*$.}\smallskip
\end{Cor}

Sous la conjecture d'Iwasawa, qui postule la trivialité de l'invariant $\mu =\mu^\emptyset_\emptyset$ ou, ce qui revient au même d'après les formules précédentes, celle de l'invariant $\mu_L$, donc en particulier lorsque le corps $K$ est absolument abélien en vertu d'un résultat célèbre de Ferrero et Washington (cf. eg. \cite{wa}), il est facile de voir que tous les paramètres $\mu_S^T$ sont en fait identiquement nuls~:

\begin{Sco} Pour tous $S$ et $T$, les paramètres $\mu_S^T$ satisfont l'inégalité~:\smallskip

\centerline{$\mu^{T *}_S + \mu^T_S \leq \mu^L + \mu_L = \mu_L^* + \mu_L$~;}\smallskip

\noindent En particulier, sous la conjecture d'Iwasawa $\mu_L = 0$, on a identiquement $\mu_S^T =0$.
\end{Sco}

\Preuve Rappelons que nous avons convenu de désigner par $H_{S_n}^{T_n}$ la $\l$-extension abélienne maximale $S$-décomposée $T$-ramifiée et d'exposant $\l^{n+1}$ du corps $K_n$. Notons donc $L_n$ l'ensemble des places $\l$-adiques de $K_n$, puis $T^\l_n$ le sous-ensemble formé de celles au-dessus de $T$ et $\ov{T}_{\!n}^\l = L_n \setminus T_n^\l$ son complémentaire~; et considérons le schéma d'extensions~: 
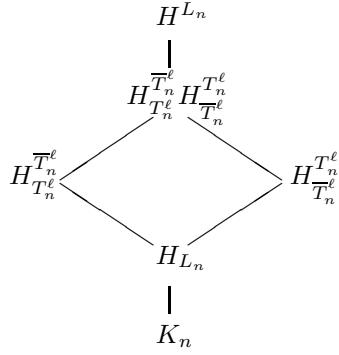
\bigskip

\begin{picture}(100,120)(0,0)
\put(155,120){$H^{L_n}$}
\put(144,90){$H^{\ov {T}_{\!n}^\l}_{T_n^\l}H^{T_n^\l}_{\ov {T}_{\!n}^\l}$}
\put(100,60){$H^{\ov{T}_{\!n}^\l}_{T_n^\l}$}
\put(205,60){$H^{T_n^\l}_{\ov {T}_{\!n}^\l}$}
\put(155,30){$H_{L_n}$}
\put(155,0){$K_n$}
\thinlines 
\put(160,114){\line(0,-1){10}}
\put(160,21){\line(0,-1){10}}
\put(154,85){\line(-3,-2){35}}
\put(167,85){\line(3,-2){35}}
\put(154,37){\line(-3,2){35}}
\put(167,37){\line(3,2){35}}
\end{picture}
\bigskip

Dans la suite exacte de $\Zl[\D]$-modules donnée par la théorie de Galois~:\smallskip

\centerline{$\Gal (H^{L_n}/K_n) \rightarrow \Gal (H^{\ov{T}_{\!n}^\l}_{T_n^\l}/K_n) \times \Gal (H_{\ov{T}_{\!n}^\l}^{T_n^\l}/K_n) \rightarrow \Gal (H_{L_n}/K_n) \rightarrow 1$,}\smallskip

\noindent tous sont paramétrés~; et il va donc de même du noyau $\Gal (H^{L_n}/ H^{\ov {T}_{\!n}^\l}_{T_n^\l}H^{T_n^\l}_{\ov {T}_{\!n}^\l})$.\par

Le Théorème 2.7 montre que caractère $\rho$ correspondant est nul. Le paramètre dominant est ainsi le caractère $\mu$ qui est donc positif ou nul. Et il vient donc~:\smallskip

\centerline{$\mu^{\ov{T}^\l}_{T^\l} + \mu^{T^\l}_{\ov{T}^\l}\ \leq \ \mu^L + \mu_L$~;}\smallskip

\noindent ce qui est précisément le résultat annoncé, en vertu du Théorème 2.8.\medskip

Venons en maintenant aux invariants $\lambda$. Le résultat principal de \cite{jm} se généralise comme suit~:

\begin{Th} Pour $S$ et $T$ arbitraires le paramètre $\lambda^T_S$ du $\l$-groupe des $S$-classes $T$-infinitésimales est donné  à partir de sa partie sauvage par la formule~:\smallskip

\centerline{$\lambda_S^T = \lambda_{S^\ell}^{T^\ell} + \chi_{T^0}^*$.}\smallskip

\noindent Il s'ensuit que la ramification modérée est déployée dans l'extension $H_{S_\infty}^{T_\infty}/K_\infty$.
\end{Th}

\Preuve Quitte à grossir l'ensemble $S$, nous pouvons supposer dans un premier temps l'inclusion $S \cup T \supset L$. Dans ce cas le résultat annoncé résulte directement des identités du miroir et du fait évident que les places modérées de $S$ sont sans influence sur l'invariant $\lambda^T_S$, puisque la montée dans la tour cyclotomique a épuisé toute possibilité d'inertie aux places modérées. Nous avons en effet~:\smallskip

\centerline{$\lambda_S^T  = \lambda_{S^\l}^T  = [\lambda_T^{S^\l} + \chi_T -1]^*- (\chi _{S^\l} - 1 ) = [\lambda_{T^\l}^{S^\l} + \chi_{T^\l} -1 ]^* + \chi_{T^0}^* - (\chi _{S^\l} - 1 )$,}\smallskip

\noindent d'où, le résultat attendu, en appliquant une seconde fois l'identité du miroir. Il s'ensuit que la ramification modérée est bien déployée (i.e. que les sous-groupes de ramification modérée sont en somme directe) dans le groupe $\Gal (H_{S_\infty}^{T_\infty}/K_\infty)$, dès que $S$ est assez grand (pour qu'on ait $S \cup T \supset L$), donc {\it a fortiori} pour $S$ petit, de sorte qu'en fin de compte la formule obtenue vaut quel que soit $S$.\medskip

Pour aller plus loin, distinguons composantes réelles et imaginaires~:\smallskip

\begin{Prop} Soient $S$ et $T$ arbitraires comme plus haut. Alors~:\smallskip

(i) La partie imaginaire du paramètre $\lambda_{S^\ell}^{T^\ell}$ est donnée par la formule~:\smallskip

\centerline{$\lambda_{S^\ell}^{T^\ell \,\ominus} = \lambda^{T^\ell \,\ominus} -  \chi_{S^\l}^\ominus$.}\smallskip

\noindent En particulier, la contribution imaginaire de la décomposition attachée aux places sauvages  non ramifiées est déployée dans l'extension $H_{S_\infty}^{T_\infty}/K_\infty$.\smallskip

(ii) Et, sous la conjecture de Leopoldt, la partie réelle de $\lambda_{S_\ell}^{T_\ell}$ est donnée par~:\smallskip

\centerline{$\lambda_{S^\ell}^{T^\ell \,\oplus} = \lambda^{T^\ell \,\oplus}$.}\smallskip

\noindent En d'autres termes, la décomposition attachée aux places sauvages non ramifiées est sans incidence sur la composante réelle de l'invariant lambda\footnote{Ce résultat est en parfaite cohérence avec la conjecture de Greenberg qui postule la nullité de l'invariant $\lambda^\oplus$ donc en particulier l'égalité $\lambda^\oplus =\lambda_{S^\l}^\oplus =0$.}.\smallskip

(iii) En fin de compte, sous la conjecture de Leopoldt, il vient ainsi~:\smallskip

\centerline{$\lambda_{S^\ell}^{T^\ell} = \lambda^{T^\ell} -  \chi_{S^\l}^\ominus$.}\smallskip

\end{Prop}

Avant d'établir cette proposition à l'aide de la Théorie des Genres dans l'appendice ci-après, tirons en tout de suite l'expression du caractère $\lambda^T_S$~:\smallskip

\begin{Th} Pour $S$ et $T$ arbitraires le paramètre $\lambda_S^T$ du $\l$-groupe des $S$-classes $T$-infinitésimales est donné  sous la conjecture de Leopoldt par la formule~:\smallskip

\centerline{$\lambda_S^T = \lambda^{T^\l} + \chi_{T^0}^* - \chi_{S^\l}^\ominus$.}\smallskip
\end{Th}

Tout comme pour les invariants $\mu^T_S$, les caractères $\lambda^T_S$ sont ainsi donnés à partir des seuls caractères $\lambda^{T^\l}$, attachés aux parties finies $T^\l$ de $L$, par des formules explicites ne faisant intervenir que des paramètres galoisiens.\par

Et, ici encore, il est possible d'exprimer les $\lambda^T_S$ à l'aide des seuls invariants référents réels. Par les identités du miroir, on a en effet~:

\begin{Sco} Lorsque $L = \ov{T}^\l \cup T^\l$ est une partition de l'ensemble $L$ des places $\l$-adiques, les paramètres lambda respectivement attachés à $\ov{T}^\l$ et à $T^\l$ satisfont les identités de dualité~:

\centerline{$\lambda^{\ov{T}^\l \ominus} = [\lambda^{T^\l \oplus} + (\chi_{T^\l}^{\,\oplus}-1)]^*$.}
\end{Sco}

\Preuve Cela résulte immédiatement du Théorème 2.6 et du Théorème 2.12.\medskip

En fin de compte, tout comme dans le cas des paramètres $\mu^T_S$, on voit que tous les paramètres $\lambda^T_S$ sont explicitement connus dès que le sont les seuls paramètres réels $\lambda^{T^\l \oplus}$ attachés aux sous-ensembles $T^\l$ de l'ensemble des places $\l$-adiques du corps considéré~:

\begin{Cor} Avec les notations du Théorème et sous les mêmes hypothèses, il vient~:

\centerline{$\lambda_S^T = \lambda^{T^\l \oplus} + [\lambda^{\ov{T}^\l \oplus} + (\chi_{\ov{T}^\l}^{\,\oplus}-1)]^* + \chi_{T^0}^* - \chi_{S^\l}^\ominus$.}\smallskip

\end{Cor}

\setcounter{section}{3}  \setcounter{Th}{0} 

\bigskip

\noindent{\large \bf 3. Appendice : généralisation d'un théorème de Greenberg}

\bigskip


Nous nous proposons ici d'établir la Proposition 2.12 et de montrer que ce résultat provient essentiellement de la Théorie des Genres. Expliquons pourquoi en examinant séparément le cas des composantes imaginaires et celui des composantes réelles~:\medskip

$(i)$ Le cas imaginaire est techniquement le plus simple. Fixons $n_0$ assez grand pour que les places $\l$-adiques se ramifient totalement dans la tour $K_\infty/K$ au-dessus de $K_{n_0}$ et convenons de noter $\mathcal L_{S_n^\l}$ le sous-groupe ambige\footnote{i.e. invariant par le groupe de Galois $\Gamma_{n_0}=\Gal (K_\infty/K_{n_0})$.} du $\l$-groupe $\Cl_{K_n}$ des classes d'idéaux (au sens ordinaire) du corps $K_n$ qui est engendré par les classes des idéaux construits sur les places au-dessus de $S^\l$ (de sorte qu'on a~: $\Cl_{S_n^\l} = \Cl_{K_n}/ \mathcal L_{S^\l_n}$).
Observons maintenant qu'il n'y a pas de capitulation pour les classes imaginaires dans la tour $K_\infty/K_{n_0}$, puisque les unités imaginaires se réduisent aux racines de l'unité, lesquelles sont cohomologiquement triviales~; de sorte que, le morphisme d'extension $\ \Cl_{K_n}^{\,\ominus}\to \Cl_{K_\infty}^{\,\ominus}$ étant injectif, la réunion $\mathcal L_{S^\l_\infty}^{\,\ominus}$ des $\mathcal L_{S^\l_n}^{\,\ominus}$ s'identifie comme $\Zl[\D]$-module au produit $(\Ql / \Zl) \otimes_{\Zl} \Zl[\D]^\ominus$. En d'autres termes, la limite projective des $\mathcal L_{S^\l_n}^{\,\ominus}$ est un $\Zl[\D]$-module de caractère $\chi_{S^\l}^\ominus$ et nous avons bien, comme attendu, l'identité~:\smallskip

\centerline{$\lambda_{S^\l}^\ominus = \lambda^\ominus - \chi_{S^\l}^\ominus$.}\smallskip

\noindent Il suit de là que la décomposition aux places sauvages (non ramifiées) est bien déployée dans la composante imaginaire du groupe de Galois $\C_{K_\infty} = \varprojlim \Cl_{K_n}$ attaché à la $\l$-extension abélienne non ramifiée maximale de $K_\infty$ ~; résultat qui est {\it a fortiori} vrai dans le groupe de Galois $\C_{K_\infty}^{T^\l}$ qui correspond, lui, à celle qui est $T^\l$-ramifiée. En fin de compte, nous avons donc tout aussi bien~:\smallskip

\centerline{$\lambda_{S_\ell}^{T_\ell \,\ominus} = \lambda^{T_\ell \,\ominus} -  \chi_{S_\l}^\ominus$,}\smallskip

\noindent ce qui est précisément le premier résultat annoncé.\medskip

(ii) Venons en maintenant au cas réel, qui est plus compliqué car il met en jeu de façon subtile la conjecture de Leopoldt. Il s'agit de voir que dans ce dernier contexte non seulement il n'y a plus injectivité du morphisme d'extension $\ \Cl_{K_n}^{\,\oplus} \to \Cl_{K_\infty}^{\,\oplus}$, mais que, de fait, la limite inductive $\mathcal L_{S^\l_\infty}^{\,\oplus}$ des $\mathcal L_{S^\l_n}^{\,\oplus}$ est nulle, de sorte que l'on a~:

\centerline{$\lambda_{S^\l}^\oplus = \lambda^\oplus$.}\smallskip 

Mieux encore, nous allons voir que la limite inductive $\mathcal L_{S^\l_\infty}^{T_\infty^\l \oplus} = \varinjlim \mathcal L_{S^\l_n}^{T_n^\l \oplus}$ des composantes réelles des sous-groupes des $\l$-groupes de classes $T_n^\l$-infinitési\-males  des corps $K_n$ respectivement engendrés par les classes des places au-dessus de $S^\l$ est encore nulle~; et que, par suite, la limite projective correspondante $ \varprojlim \mathcal L_{S^\l_n}^{T_n^\l \oplus}$ est pseudo-nulle, de sorte que nous avons bien, comme annoncé~:\medskip

\centerline{$\lambda_{S_\ell}^{T_\ell \,\oplus} = \lambda^{T_\ell \,\oplus}$.}\bigskip

Enonçons explicitement ce dernier résultat, qui généralise dans un contexte $\l$-ramifié un théorème de semi-simplicité dû à R. Greenberg (cf. \cite{gb})~:\smallskip

\begin{Th} Soit $K^+$ un corps de nombres totalement réel satisfaisant la conjecture de Leopoldt pour un premier $\l$ (i.e. dont le rang $\l$-adique du goupe des unités est égal à $[K^+:\Q\,]-1$) et $K_\infty^+ = \cup_{n\in\N} K_n^+$ sa $\Zl$-extension cyclotomique. Alors, pour toute partition $L = S\cup T$ de l'ensemble $L$ des places $\l$-adiques de $K^+$, le sous-groupe $\mathcal L_{S_\infty}^{T_\infty}$  du $\l$-groupe $\Cl^{\,T_\infty}$ des classes $T$-infinitésimales du corps $K^+_\infty$ qui est engendré par les classes des idéaux au-dessus de $S$ est trivial, dès lors que les places de $S$ ne se décomposent pas dans l'extension $K^+_\infty/K^+$.
\end{Th}

\Remarque La condition finale de non décomposition étant banalement vérifiée à partir d'un étage fini $K^+_{n_0}$ de la tour cyclotomique, la conclusion du théorème reste donc valide dès lors que le corps $K^+_{n_0}$ satisfait la conjecture de Leopoldt.\medskip

\Preuve Comme nous l'avons dit, la démonstration de ce résultat repose sur la Théorie des Genres~: l'hypothèse de non décomposition nous dit que les idéaux construits sur les places de $S$ sont ambiges dans l'extension $K^+_\infty/K^+$ (i.e. invariantes par l'action de $\Gamma = \Gal(K^+_\infty/K^+)$), donc que leurs images respectives dans chacun des $\l$-groupes de classes $T$-infinitésimales $\Cl_{K^+_n}^T$ tombent dans le sous-groupe ambige $(\Cl_{K^+_n}^T)^\Gamma$. Tout revient donc à établir que ces sous-groupes ont un ordre borné indépendamment de $n$ sous les hypothèses du Théorème ou, ce qui revient au même, que les quotients des genres ${}^\Gamma(\Cl_{K^+_n}^T)$ sont d'ordre borné.\par

Observons au passage que si la conjecture de Leopoldt se trouve vérifiée dans la tour $K_\infty^+$, i.e. si {\it tous} les corps $K^+_n$ la satisfont, les groupes $\Cl_{K^+_n}^T$ de $\l$-classes $T$-infinitésimales sont alors finis~; de sorte que dans ce cas le sous-groupe ambige $(\Cl_{K^+_n}^T)^\Gamma$ a même ordre que le quotient des genres ${}^\Gamma(\Cl_{K^+_n}^T)$ en vertu de la suite exacte canonique associée à $\omega_0 = \gamma-1$~:
$$\begin{CD}
1 @>>> (\Cl_{K^+_n}^T)^\Gamma @>>> \Cl_{K^+_n}^T @>\omega_0>>  \Cl_{K^+_n}^T @>>> {}^\Gamma(\Cl_{K^+_n}^T) @>>> 1
\end{CD}$$
\noindent mais cette hypothèse est inutile à notre conclusion, le seul point qui importe étant que les sous-groupes ambiges et les quotients des genres sont simultanément bornés ou non\footnote{Ce qui se lit sur le polynôme caractéristique de la limite projective des $\,\Cl_{K^+_n}^T$.}. Cela étant, le Théorème annoncé résulte de l'expression du nombre de genres donnée par la suite exacte des genres $T$-infinitésimaux que nous allons maintenant énoncer~:

\begin{Th} Soit $F$ un corps de nombres, $F_\infty = \cup_{n\in\N}F_n$ sa $\Zl$-extension cyclotomique et $S$ une partie non vide de l'ensemble $L$ des places $\l$-adiques de $F$ de complémentaire $T=L\setminus S$. Pour chaque entier $n$, désignons par $H^T_{F_n}$ la $\l$-extension abélienne $T$-ramifiée maximale de $F_n$~; notons $H^T_{F_\infty}=\cup_{n\in\N}\,H^T_{F_n}$ leur réunion et $H^{T\,\rm ab}_{F_\infty} = \cup_{n \in \N}\,H^T_{F_n /F}$ la sous-extension maximale de  $H^T_{F_\infty}$ qui est abélienne sur $F$. Avec ces notations, le groupe de Galois $\G^T_{F_\infty /F} = \Gal(H^{T\,\rm ab}_{F_\infty}/F_\infty)$ est donné par la suite exacte des $T$-genres infinitésimaux~:
$$
\E^T_F / (\E^T_F \cap \mathcal N_{F_\infty / F}) \hookrightarrow \wi \oplus_{\p\in S} I_\p(H^{T\,\rm ab}_{F_\infty} /F) \to \G^T_{F_\infty /F}  \to \Cl^T_F \twoheadrightarrow \Gal (H^T_F \cap F_\infty /F)
$$
Dans celle-ci $\,\E^T_F$ est le $\l$-groupe des unités $T$-infinitésimales du corps $F$ et $\mathcal N_{F_\infty / F}$ le groupe des normes cyclotomiques~; la somme au centre porte sur les familles $(\sigma_\p)$ d'éléments des sous-groupes procycliques d'inertie attachés aux places de $S$ qui satisfont la formule du produit $\prod_{\p\in S}\sigma_\p \mid_{F_\infty} =1$~; enfin $\,\Cl^T_F \simeq \Gal(H^T_F/F)$ est le $\l$-groupe des classes $T$-infinitésimales du corps $F$.\smallskip
\end{Th}

\begin{Cor} Lorsque $F$ est totalement réel et vérifie la conjecture de Leopoldt, le quotient $\G^T_{F_\infty /F}$ des  $T$-genres infinitésimaux de la tour cyclotomique $F_\infty/F$ est fini. En d'autres termes, les quotients des genres $\G^T_{F_n /F}  = {}^\Gamma \Cl^T_{F_n}$ attachés aux étages finis de la tour sont eux-mêmes finis et ultimement constants~; le nombre de classes ambiges $T$-infinitésimales $\mid \Cl^{T\ \Gamma}_{F_n} \mid$est stationnaire dans la tour~; et par suite, dès lors que les places de $S$ ne se décomposent pas dans le premier étage de la tour, les classes des idéaux au-dessus de $S$ capitulent dans le groupe $\,\Cl^T_{F_\infty}$.  
\end{Cor}

\Preuve Commençons par établir l'exactitude de la suite des genres. Le conoyau de la restriction canonique $\Gal (H^{T\,\rm ab}_{F_\infty}/F_\infty) \to \Gal (H^T_F/F)$ est évidemment le groupe $\Gal (H^T_F \cap F_\infty /F)$, qui est fini, puisque les places de $S$ se ramifient dans la tour cyclotomique $F_\infty/F$. Et son noyau est bien engendré par les sous-groupes d'inertie des places de $S$, lesquels sont procycliques puisque ces places sont presque totalement ramifiées dans la tour cyclotomique.\par
Reste à voir que le noyau à gauche est l'image du groupe $\,\E_F^T$ des unités $T$-infinitésimales du corps $F$, résultat qui peut être regardé comme la transposition dans notre contexte du Théorème 6 de \cite{j-2}~: le point essentiel consiste à voir que toute famille $(\sigma_\p)_\p$ qui vérifie la formule du produit provient bien, par les symboles de restes normiques de Hasse, d'une unité $T$-infinitésimale. Or les arguments de \cite{j-2} appliqués aux extensions abéliennes finies $F_n/F$ montrent, après restriction de la condition d'infinitésimalité aux places de $L$ contenues dans $T$, que c'est bien le cas à chaque niveau fini de la tour (i.e. dans $\wi \oplus_{\p\in S} I_\p(H^T_{F_n /F})$~; de sorte que  l'on conclut aisément par compacité en passant à la limite projective.\smallskip

Cela acquis, le Corollaire en découle sans difficulté~: écrivons $t = \sum_{\p \in T} [F_\p : \Q_\l]$ le degré total en $T$ et $s = \sum_{\p \in S} [F_\p : \Q_\l]$ celui en $S$, de sorte que nous avons $t+s=r=[F :\Q]$. Notons ensuite $\wi \U_L$ le noyau\footnote{On prendra garde que {\it ce n'est pas} le le groupe des classes logarithmiques au sens de \cite{ja2}.} dans $\U_L = \prod_{\p \in L}\U_\p$ de la norme arithmétique $N_{F/ \Q}$ et écrivons de même $\wi \U_L$ le noyau de la norme $N_{F/ \Q}$ dans $\U_S = \prod_{\p \in S}\U_\p$. Sous la conjecture de Leopoldt dans le corps $F$, le $\l$-adifié $\,\E_F = \Zl \otimes E_F$ du groupe des unités de $F$ (au sens habituel) est un $\Zl$-module de rang $r-1$ qui s'injecte dans $\,\wi \U_L$ avec un indice fini. Il en résulte que son sous-groupe $T$-infinitésimal $\,\E^T_F$, qui est le noyau du morphisme naturel de $\,\E_F$ dans $\,\U_T = \prod_{\p \in T} \U_\p$, s'injecte dans $\,\wi \U_S$ avec un indice fini. D'où le résultat attendu, la somme $ \wi \oplus_{\p\in S} I_\p(H^T_{F_\infty /F})$ étant précisément l'image canonique de $\,\wi \U_S$ par les symboles locaux de réciprocité.


{\small
\def\refname{
 }

\bigskip\noindent
{\small
\begin{tabular}{l}
{Jean-François {\sc Jaulent}}\\
Institut de Mathématiques de Bordeaux\\
Université Bordeaux I \\
351, cours de la Libération\\
F-33405 Talence Cedex\\
email : jaulent@math.u-bordeaux.fr 
\end{tabular}
}

\end{document}